
\input amstex



\def\spaces{\space\space\space\space\space\space\space\space\space\space}
\def\spacess{\message{\spaces\spaces\spaces\spaces\spaces\spaces\spaces}}
\spacess
\spacess
\message{Annals of Mathematics Style: Current Version: 1.1. June 10, 1992}
\spacess
\spacess

\catcode`\@=11

\hyphenation{acad-e-my acad-e-mies af-ter-thought anom-aly anom-alies
an-ti-deriv-a-tive an-tin-o-my an-tin-o-mies apoth-e-o-ses apoth-e-o-sis
ap-pen-dix ar-che-typ-al as-sign-a-ble as-sist-ant-ship as-ymp-tot-ic
asyn-chro-nous at-trib-uted at-trib-ut-able bank-rupt bank-rupt-cy
bi-dif-fer-en-tial blue-print busier busiest cat-a-stroph-ic
cat-a-stroph-i-cally con-gress cross-hatched data-base de-fin-i-tive
de-riv-a-tive dis-trib-ute dri-ver dri-vers eco-nom-ics econ-o-mist
elit-ist equi-vari-ant ex-quis-ite ex-tra-or-di-nary flow-chart
for-mi-da-ble forth-right friv-o-lous ge-o-des-ic ge-o-det-ic geo-met-ric
griev-ance griev-ous griev-ous-ly hexa-dec-i-mal ho-lo-no-my ho-mo-thetic
ideals idio-syn-crasy in-fin-ite-ly in-fin-i-tes-i-mal ir-rev-o-ca-ble
key-stroke lam-en-ta-ble light-weight mal-a-prop-ism man-u-script
mar-gin-al meta-bol-ic me-tab-o-lism meta-lan-guage me-trop-o-lis
met-ro-pol-i-tan mi-nut-est mol-e-cule mono-chrome mono-pole mo-nop-oly
mono-spline mo-not-o-nous mul-ti-fac-eted mul-ti-plic-able non-euclid-ean
non-iso-mor-phic non-smooth par-a-digm par-a-bol-ic pa-rab-o-loid
pa-ram-e-trize para-mount pen-ta-gon phe-nom-e-non post-script pre-am-ble
pro-ce-dur-al pro-hib-i-tive pro-hib-i-tive-ly pseu-do-dif-fer-en-tial
pseu-do-fi-nite pseu-do-nym qua-drat-ics quad-ra-ture qua-si-smooth
qua-si-sta-tion-ary qua-si-tri-an-gu-lar quin-tes-sence quin-tes-sen-tial
re-arrange-ment rec-tan-gle ret-ri-bu-tion retro-fit retro-fit-ted
right-eous right-eous-ness ro-bot ro-bot-ics sched-ul-ing se-mes-ter
semi-def-i-nite semi-ho-mo-thet-ic set-up se-vere-ly side-step sov-er-eign
spe-cious spher-oid spher-oid-al star-tling star-tling-ly
sta-tis-tics sto-chas-tic straight-est strange-ness strat-a-gem strong-hold
sum-ma-ble symp-to-matic syn-chro-nous topo-graph-i-cal tra-vers-a-ble
tra-ver-sal tra-ver-sals treach-ery turn-around un-at-tached un-err-ing-ly
white-space wide-spread wing-spread wretch-ed wretch-ed-ly Brown-ian
Eng-lish Euler-ian Feb-ru-ary Gauss-ian Grothen-dieck Hamil-ton-ian
Her-mit-ian Jan-u-ary Japan-ese Kor-te-weg Le-gendre Lip-schitz
Lip-schitz-ian Mar-kov-ian Noe-ther-ian No-vem-ber Rie-mann-ian
Schwarz-schild Sep-tem-ber Za-mo-lod-chi-kov Knizh-nik quan-tum Op-dam
Mac-do-nald Ca-lo-ge-ro Su-ther-land Mo-ser Ol-sha-net-sky  Pe-re-lo-mov
in-de-pen-dent ope-ra-tors Hec-man Op-dam har-mo-nic ana-ly-sis
al-geb-ra pre-print char-ac-ter-is-tic sub-mo-dule re-la-tion Heck-man
Ha-rish-Chan-dra Hel-ga-son Su-ther-land em-bed-ding
}

\Invalid@\nofrills
\Invalid@\usualspace
\newif\ifnofrills@
\def\nofrills@#1#2{\relaxnext@
  \DN@{\ifx\next\nofrills
    \nofrills@true\let#2\relax\DN@\nofrills{\nextii@}%
  \else
    \nofrills@false\def#2{#1}\let\next@\nextii@\fi
\next@}}
\def\usualspace@#1{\ifnofrills@\def\usualspace{#1}\fi}
\def\addto#1#2{\csname \expandafter\eat@\string#1@\endcsname
  \expandafter{\the\csname \expandafter\eat@\string#1@\endcsname#2}}
\newdimen\bigsize@
\def\big@#1#2{{\hbox{$\left#2\vcenter to#1\bigsize@{}%
  \right.\nulldelimiterspace\z@\m@th$}}}
\def\big{\big@\@ne}
\def\Big{\big@{1.5}}
\def\bigg{\big@\tw@}
\def\Bigg{\big@{2.5}}
\def\raggedcenter@{\leftskip\z@ plus.4\hsize \rightskip\leftskip
 \parfillskip\z@ \parindent\z@ \spaceskip.3333em \xspaceskip.5em
 \pretolerance9999\tolerance9999 \exhyphenpenalty\@M
 \hyphenpenalty\@M \let\\\linebreak}
\def\upperspecialchars{\def\ss{SS}\let\i=I\let\j=J\let\ae\AE\let\oe\OE
  \let\o\O\let\aa\AA\let\l\L}
\def\uppercasetext@#1{%
  {\spaceskip1.2\fontdimen2\the\font plus1.2\fontdimen3\the\font
   \upperspecialchars\uctext@#1$\m@th\aftergroup\eat@$}}
\def\uctext@#1$#2${\endash@#1-\endash@$#2$\uctext@}
\def\endash@#1-#2\endash@{%
\uppercase{#1}\if\notempty{#2}--\endash@#2\endash@\fi}
\def\runaway@#1{\DN@{#1}\ifx\envir@\next@
  \Err@{You seem to have a missing or misspelled \string\end#1 ...}%
  \let\envir@\empty\fi}
\newif\iftemp@
\def\notempty#1{TT\fi\def\test@{#1}\ifx\test@\empty\temp@false
  \else\temp@true\fi \iftemp@}

\font@\tensmc=cmcsc10
\font@\sevenex=cmex7 
\font@\sevenit=cmti7
\font@\eightrm=cmr8 
\font@\sixrm=cmr6 
\font@\eighti=cmmi8     \skewchar\eighti='177 
\font@\sixi=cmmi6       \skewchar\sixi='177   
\font@\eightsy=cmsy8    \skewchar\eightsy='60 
\font@\sixsy=cmsy6      \skewchar\sixsy='60   
\font@\eightex=cmex8 %
\font@\eightbf=cmbx8 
\font@\sixbf=cmbx6   
\font@\eightit=cmti8 
\font@\eightsl=cmsl8 
\font@\eightsmc=cmcsc10 
\font@\eighttt=cmtt8 

\loadmsam
\loadmsbm
\loadeufm
\UseAMSsymbols

\def\penaltyandskip@#1#2{\relax\ifdim\lastskip<#2\relax\removelastskip
      \ifnum#1=\z@\else\penalty@#1\relax\fi\vskip#2%
  \else\ifnum#1=\z@\else\penalty@#1\relax\fi\fi}
\def\nobreak{\penalty\@M
  \ifvmode\def\penalty@{\let\penalty@\penalty\count@@@}%
  \everypar{\let\penalty@\penalty\everypar{}}\fi}
\let\penalty@\penalty

\def\block{\RIfMIfI@\nondmatherr@\block\fi
       \else\ifvmode\vskip\abovedisplayskip\noindent\fi
        $$\def\endblock{\par\egroup$$}\fi
  \vbox\bgroup\advance\hsize-2\indenti\noindent}
\def\endblock{\par\egroup}

\def\logo@{\baselineskip2pc \hbox to\hsize{\hfil\eightpoint Typeset by
 \AmSTeX}}




\font\elevensc=cmcsc10 scaled\magstephalf
\font\tensc=cmcsc10

\font\eightsc=cmcsc10 scaled800

\font\elevenrm=cmr10 scaled \magstephalf
\font\ninerm=cmr9
\font\eightrm=cmr8
\font\sixrm=cmr6
\font\fiverm=cmr5

\font\eleveni=cmmi10 scaled\magstephalf
\font\ninei=cmmi9
\font\eighti=cmmi8
\font\sixi=cmmi6
\font\fivei=cmmi5
\skewchar\ninei='177 \skewchar\eighti='177 \skewchar\sixi='177
\skewchar\eleveni='177

\font\elevensy=cmsy10 scaled\magstephalf
\font\ninesy=cmsy9
\font\eightsy=cmsy8
\font\sixsy=cmsy6
\font\fivesy=cmsy5
\skewchar\ninesy='60 \skewchar\eightsy='60 \skewchar\sixsy='60
\skewchar\elevensy'60

\font\eighteenbf=cmbx10 scaled\magstep3

\font\twelvebf=cmbx10 scaled \magstep1
\font\elevenbf=cmbx10 scaled \magstephalf
\font\tenbf=cmbx10
\font\ninebf=cmbx9
\font\eightbf=cmbx8
\font\sixbf=cmbx6
\font\fivebf=cmbx5

\font\elevenit=cmti10 scaled\magstephalf
\font\nineit=cmti9
\font\eightit=cmti8

\font\eighteenmib=cmmib10 scaled \magstep3
\font\twelvemib=cmmib10 scaled \magstep1
\font\elevenmib=cmmib10 scaled\magstephalf
\font\tenmib=cmmib10
\font\eightmib=cmmib10 scaled 800 
\font\sixmib=cmmib10 scaled 600

\font\eighteensyb=cmbsy10 scaled \magstep3 
\font\twelvesyb=cmbsy10 scaled \magstep1
\font\elevensyb=cmbsy10 scaled \magstephalf
\font\tensyb=cmbsy10 
\font\eightsyb=cmbsy10 scaled 800
\font\sixsyb=cmbsy10 scaled 600
 
\font\elevenex=cmex10 scaled \magstephalf
\font\tenex=cmex10     
\font\eighteenex=cmex10 scaled \magstep3


\def\elevenpoint{\def\rm{\fam0\elevenrm}%
  \textfont0=\elevenrm \scriptfont0=\eightrm \scriptscriptfont0=\sixrm
  \textfont1=\eleveni \scriptfont1=\eighti \scriptscriptfont1=\sixi
  \textfont2=\elevensy \scriptfont2=\eightsy \scriptscriptfont2=\sixsy
  \textfont3=\elevenex \scriptfont3=\tenex \scriptscriptfont3=\tenex
  \def\bf{\fam\bffam\elevenbf}%
  \def\it{\fam\itfam\elevenit}%
  \textfont\bffam=\elevenbf \scriptfont\bffam=\eightbf
   \scriptscriptfont\bffam=\sixbf
\normalbaselineskip=13.95pt
  \setbox\strutbox=\hbox{\vrule height9.5pt depth4.4pt width0pt\relax}%
  \normalbaselines\rm}

\elevenpoint 

\def\ninepoint{\def\rm{\fam0\ninerm}%
  \textfont0=\ninerm \scriptfont0=\sixrm \scriptscriptfont0=\fiverm
  \textfont1=\ninei \scriptfont1=\sixi \scriptscriptfont1=\fivei
  \textfont2=\ninesy \scriptfont2=\sixsy \scriptscriptfont2=\fivesy
  \textfont3=\tenex \scriptfont3=\tenex \scriptscriptfont3=\tenex
  \def\it{\fam\itfam\nineit}%
  \textfont\itfam=\nineit
  \def\bf{\fam\bffam\ninebf}%
  \textfont\bffam=\ninebf \scriptfont\bffam=\sixbf
   \scriptscriptfont\bffam=\fivebf
\normalbaselineskip=11pt
  \setbox\strutbox=\hbox{\vrule height8pt depth3pt width0pt\relax}%
  \normalbaselines\rm}

\def\eightpoint{\def\rm{\fam0\eightrm}%
  \textfont0=\eightrm \scriptfont0=\sixrm \scriptscriptfont0=\fiverm
  \textfont1=\eighti \scriptfont1=\sixi \scriptscriptfont1=\fivei
  \textfont2=\eightsy \scriptfont2=\sixsy \scriptscriptfont2=\fivesy
  \textfont3=\tenex \scriptfont3=\tenex \scriptscriptfont3=\tenex
  \def\it{\fam\itfam\eightit}%
  \textfont\itfam=\eightit
  \def\bf{\fam\bffam\eightbf}%
  \textfont\bffam=\eightbf \scriptfont\bffam=\sixbf
   \scriptscriptfont\bffam=\fivebf
\normalbaselineskip=12pt
  \setbox\strutbox=\hbox{\vrule height8.5pt depth3.5pt width0pt\relax}%
  \normalbaselines\rm}


\def\eighteenbold{\def\rm{\fam0\eighteenbf}%
  \textfont0=\eighteenbf \scriptfont0=\twelvebf \scriptscriptfont0=\tenbf
  \textfont1=\eighteenmib \scriptfont1=\twelvemib\scriptscriptfont1=\tenmib
  \textfont2=\eighteensyb \scriptfont2=\twelvesyb\scriptscriptfont2=\tensyb
  \textfont3=\eighteenex \scriptfont3=\tenex \scriptscriptfont3=\tenex
  \def\bf{\fam\bffam\eighteenbf}%
  \textfont\bffam=\eighteenbf \scriptfont\bffam=\twelvebf
   \scriptscriptfont\bffam=\tenbf
\normalbaselineskip=20pt
  \setbox\strutbox=\hbox{\vrule height13.5pt depth6.5pt width0pt\relax}%
\everymath {\fam0 }
\everydisplay {\fam0 }
  \normalbaselines\rm}

\def\elevenbold{\def\rm{\fam0\elevenbf}%
  \textfont0=\elevenbf \scriptfont0=\eightbf \scriptscriptfont0=\sixbf
  \textfont1=\elevenmib \scriptfont1=\eightmib \scriptscriptfont1=\sixmib
  \textfont2=\elevensyb \scriptfont2=\eightsyb \scriptscriptfont2=\sixsyb
  \textfont3=\elevenex \scriptfont3=\elevenex \scriptscriptfont3=\elevenex
  \def\bf{\fam\bffam\elevenbf}%
  \textfont\bffam=\elevenbf \scriptfont\bffam=\eightbf
   \scriptscriptfont\bffam=\sixbf
\normalbaselineskip=14pt
  \setbox\strutbox=\hbox{\vrule height10pt depth4pt width0pt\relax}%
\everymath {\fam0 }
\everydisplay {\fam0 }
  \normalbaselines\bf}

\hsize=31pc
\vsize=48pc

\parindent=22pt
\parskip=0pt

\widowpenalty=10000
\clubpenalty=10000

\topskip=12pt 

\skip\footins=20pt
\dimen\footins=3in 

\abovedisplayskip=6.95pt plus3.5pt minus 3pt
\belowdisplayskip=\abovedisplayskip


\voffset=7pt\hoffset= .7in

\newif\iftitle

\def\amheadline{\iftitle%
\hbox to\hsize{\hss\currannalsline\hss}\else\line{\ifodd\pageno
\hfill\thetitle\hfill\llap{\elevenrm\folio}\else\rlap{\elevenrm\folio}
\hfill\theauthors\hfill\fi}\fi}

\headline={\amheadline}
\footline={\global\titlefalse}


\def\annalsline#1#2{\vfill\eject
\ifodd\pageno\else 
\line{\hfill}
\vfill\eject\fi
\global\titletrue
\def\currannalsline{\eightrm 
{\eightbf#1} (#2), \thepages}}

\def\titleheadline#1{\def\one{#1}\ifx\one\empty\else
\def\thetitle{{
\let\\ \relax\eightsc\uppercase{#1}}}\fi}

\newif\ifshort

\let\shorttitle\titleheadline

\def\onpages#1#2{\def\thepages{#1--#2}}

\def\thismuchskip[#1]{\vskip#1pt}
\def\ilook{\ifx\next[ \let\go\thismuchskip\else
\let\go\relax\vskip1pt\fi\go}

\def\institution#1{\def\theinstitutions{\vbox{\baselineskip10pt
\def\and{{\eightrm and }}
\def\\{\futurelet\next\ilook}\eightsc #1}}}
\let\institutions\institution

\newwrite\auxfile

\def\startingpage#1{\def\one{#1}\ifx\one\empty\global\pageno=1\else
\global\pageno=#1\fi
\theoremcount=0 \eqcount=0 \sectioncount=0 
\openin1 \jobname.aux \ifeof1 
\onpages{#1}{???}
\else\closein1 \relax\input \jobname.aux
\onpages{#1}{\lastpage}
\fi\immediate\openout\auxfile=\jobname.aux
}

\def\endarticle{\ifRefsUsed\global\RefsUsedfalse%
\else\vskip21pt\theinstitutions%
\nobreak\vskip8pt
\fi%
\write\auxfile{\string\def\string\lastpage{\the\pageno}}}

\outer\def\bye{\endarticle\par \vfill \supereject \end}

\def\document{\let\fontlist@\relax\let\alloclist@\relax
 \elevenpoint}


\newif\ifacks
\long\def\acknowledgements#1{\def\one{#1}\ifx\one\empty\else
\vskip-\baselineskip%
\global\ackstrue\footnote{\ \unskip}{*#1}\fi}

\def\title#1{\titleheadline{#1}
\vbox to80pt{\vfill
\baselineskip=18pt
\parindent=0pt
\overfullrule=0pt
\hyphenpenalty=10000
\everypar={\hskip\parfillskip\relax}
\hbadness=10000
\def\\ {\vskip1sp}
\eighteenbold#1\vskip1sp}}

\newif\ifauthor

\def\author#1{\vskip11pt
\hbox to\hsize{\hss\tenrm By \tensc#1\ifacks\global\acksfalse*\fi\hss}
\ifshort\else\xdef\theauthors{{\eightsc\uppercase{#1}}}\fi%
\vskip21pt\global\authortrue\everypar={\global\authorfalse\everypar={}}}

\def\twoauthors#1#2{\vskip11pt
\hbox to\hsize{\hss%
\tenrm By \tensc#1 {\tenrm and} #2\ifacks\global\acksfalse*\fi\hss}
\ifshort\else\xdef\theauthors{{\eightsc\uppercase{#1 and #2}}}\fi%
\vskip21pt
\global\authortrue\everypar={\global\authorfalse\everypar={}}}


\newcount\theoremcount
\newcount\sectioncount
\newcount\eqcount

\newif\ifspecialnumon

\def\eqnumber=#1 {\global\eqcount=#1 \global\advance\eqcount by-1\relax}
\def\sectionnumber=#1 {\global\sectioncount=#1 
\global\advance\sectioncount by-1\relax}
\def\proclaimnumber=#1 {\global\theoremcount=#1 
\global\advance\theoremcount by-1\relax}

\newif\ifsection
\newif\ifsubsection

\def\elevenboldmath#1{$#1$\egroup}
\def\mathbold{\hbox\bgroup\elevenbold\elevenboldmath}

\def\section#1{\global\theoremcount=0
\global\eqcount=0
\ifauthor\global\authorfalse\else%
\vskip18pt plus 18pt minus 6pt\fi%
{\parindent=0pt
\everypar={\hskip\parfillskip}
\def\\ {\vskip1sp}\elevenpoint\bf%
\ifspecialnumon\global\specialnumonfalse$\rm\spnum$%
\gdef\sectnum{$\rm\spnum$}%
\else\interlinepenalty=10000%
\global\advance\sectioncount by1\relax\the\sectioncount%
\gdef\sectnum{\the\sectioncount}%
\fi. \hskip6pt#1
\vrule width0pt depth12pt}
\hskip\parfillskip
\global\sectiontrue%
\everypar={\global\sectionfalse\global\interlinepenalty=0\everypar={}}%
\ignorespaces

}


\newif\ifspequation

\let\eqno\leqno 

\newif\ifineqalignno
\let\saveleqalignno\leqalignno                        
\def\leqalignno{\let\eqnu\Eeqnu\saveleqalignno}

\let\eqalignno\leqalignno

\def\sectandeqnum{%
\ifspecialnumon\global\specialnumonfalse
$\rm\spnum$\gdef\eqnum{{$\rm\spnum$}}\else\global\firstlettertrue
\global\advance\eqcount by1 
\ifappend\applett\else\the\sectioncount\fi.%
\the\eqcount
\xdef\eqnum{\ifappend\applett\else\the\sectioncount\fi.\the\eqcount}\fi}

\def\eqnu{\leqno{\hbox{\elevenrm\ifspequation\else(\fi\sectandeqnum
\ifspequation\global\spequationfalse\else)\fi}}}      

\def\Speqnu{\global\setbox\leqnobox=\hbox{\elevenrm
\ifspequation\else%
(\fi\sectandeqnum\ifspequation\global\spequationfalse\else)\fi}}

\def\Eeqnu{\hbox{\elevenrm
\ifspequation\else%
(\fi\sectandeqnum\ifspequation\global\spequationfalse\else)\fi}}

\newif\iffirstletter
\global\firstlettertrue
\def\eqletter#1{\global\specialnumontrue\iffirstletter\global\firstletterfalse
\global\advance\eqcount by1\fi
\gdef\spnum{\the\sectioncount.\the\eqcount#1}\eqnu}

\newbox\leqnobox
\def\outsideeqnu#1{\global\setbox\leqnobox=\hbox{#1}}

\def\eatone#1{}

\def\dosplit#1#2{\vskip-.5\abovedisplayskip
\setbox0=\hbox{$\let\eqno\outsideeqnu%
\let\eqnu\Speqnu\let\leqno\outsideeqnu#2$}%
\setbox1\vbox{\noindent\hskip\wd\leqnobox\ifdim\wd\leqnobox>0pt\hskip1em\fi%
$\displaystyle#1\mathstrut$\hskip0pt plus1fill\relax
\vskip1pt
\line{\hfill$\let\eqnu\eatone\let\leqno\eatone%
\displaystyle#2\mathstrut$\ifmathqed~~\qed\fi}}%
\copy1
\ifvoid\leqnobox
\else\dimen0=\ht1 \advance\dimen0 by\dp1
\vskip-\dimen0
\vbox to\dimen0{\vfill
\hbox{\unhbox\leqnobox}
\vfill}
\fi}

\everydisplay{\lookforbreak}

\long\def\lookforbreak #1$${\def\mathone{#1}
\expandafter\testforbreak\mathone\splitmath @}

\def\testforbreak#1\splitmath #2@{\def\mathtwo{#2}\ifx\mathtwo\empty%
#1$$%
\ifmathqed\vskip-\belowdisplayskip
\setbox0=\vbox{\let\eqno\relax\let\eqnu\relax$\displaystyle#1$}%
\vskip-\ht0\vskip-3.5pt\hbox to\hsize{\hfill\qed}
\vskip\ht0\vskip3.5pt\fi
\else$$\vskip-\belowdisplayskip
\vbox{\dosplit{#1}{\let\eqno\eatone
\let\splitmath\relax#2}}%
\nobreak\vskip.5\belowdisplayskip
\noindent\ignorespaces\fi}


\newif\ifmathqed



\newcount\linenum
\newcount\colnum

\def\spline{\omit&\multispan{\the\colnum}{\hrulefill}\cr}
\def\colcounter{\ifnum\linenum=1\global\advance\colnum by1\fi}

\def\everyline{\noalign{\global\advance\linenum by1\relax}%
\ifnum\linenum=2\spline\fi}

\def\mtable{\bgroup\offinterlineskip
\everycr={\everyline}\global\linenum=0
\halign\bgroup\vrule height 10pt depth 4pt width0pt
\hfill$##$\hfill\hskip6pt\ifnum\linenum>1
\vrule\fi&&\colcounter\hskip12pt\hfill$##$\hfill\hskip12pt\cr}

\def\endmtable{\crcr\egroup\egroup}




\def\xast{*}
\newcount\intable
\newcount\mathcol
\newcount\savemathcol
\newcount\topmathcol
\newdimen\arrayhspace
\newdimen\arrayvspace

\arrayhspace=8pt 
\arrayvspace=12pt 

\newif\ifdollaron

\def\mathalign#1{\def\arg{#1}\ifx\arg\xast%
\let\go\relax\else\let\go\mathalign%
\global\advance\mathcol by1 %
\global\advance\topmathcol by1 %
\expandafter\def\csname  mathcol\the\mathcol\endcsname{#1}%
\fi\go}

\def\arraypickapart#1]#2*{\if#1c \ifmmode\vcenter\else
\global\dollarontrue$\vcenter\fi\else%
\if#1t\vtop\else\if#1b\vbox\fi\fi\fi\bgroup%
\def\one{#2}}

\def\arraystrut{\vrule height .7\arrayvspace depth .3\arrayvspace width 0pt}

\def\array#1#2*{\def\firstarg{#1}%
\if\firstarg[ \def\two{#2} \expandafter\arraypickapart\two*\else%
\ifmmode\vcenter\else\vbox\fi\bgroup \def\one{#1#2}\fi%
\global\everycr={\noalign{\global\mathcol=\savemathcol\relax}}%
\def\\ {\cr}%
\global\advance\intable by1 %
\ifnum\intable=1 \global\mathcol=0 \savemathcol=0 %
\else \global\advance\mathcol by1 \savemathcol=\mathcol\fi%
\expandafter\mathalign\one*%
\mathcol=\savemathcol %
\halign\bgroup&\hskip.5\arrayhspace\arraystrut%
\global\advance\mathcol by1 \relax%
\expandafter\if\csname mathcol\the\mathcol\endcsname r\hfill\else%
\expandafter\if\csname mathcol\the\mathcol\endcsname c\hfill\fi\fi%
$\displaystyle##$%
\expandafter\if\csname mathcol\the\mathcol\endcsname r\else\hfill\fi\relax%
\hskip.5\arrayhspace\cr}

\def\endarray{\crcr\egroup\egroup%
\global\mathcol=\savemathcol %
\global\advance\intable by -1\relax%
\ifnum\intable=0 %
\ifdollaron\global\dollaronfalse $\fi
\loop\ifnum\topmathcol>0 %
\expandafter\def\csname  mathcol\the\topmathcol\endcsname{}%
\global\advance\topmathcol by-1 \repeat%
\global\everycr={}\fi%
}

\def\big#1{{\hbox{$\left#1\vbox to 10pt{}\right.\n@space$}}}
\def\Big#1{{\hbox{$\left#1\vbox to 13pt{}\right.\n@space$}}}
\def\bigg#1{{\hbox{$\left#1\vbox to 16pt{}\right.\n@space$}}}
\def\Bigg#1{{\hbox{$\left#1\vbox to 19pt{}\right.\n@space$}}}


\def\figcaption#1#2#3{\topinsert
\vskip4pt 
\vbox to#3{\vfill}\vskip1sp
\setbox0=\hbox{\eightsc Figure #1.\hskip12pt\eightpoint #2}
\ifdim\wd0>\hsize
\noindent\eightsc Figure #1.\hskip12pt\eightpoint #2
\else
\centerline{\eightsc Figure #1.\hskip12pt\eightpoint #2}
\fi
\vskip16pt
\endinsert}

\def\wfig#1#2#3{\topinsert
\vskip4pt 
\hbox to\hsize{\hss\vbox{\hrule height .25pt width #3
\hbox to #3{\vrule width .25pt height #2\hfill\vrule width .25pt height #2}
\hrule height.25pt}\hss}
\vskip1sp
\centerline{\eightsc Figure #1}
\vskip16pt
\endinsert}

\def\wfigcaption#1#2#3#4{\topinsert
\vskip4pt 
\hbox to\hsize{\hss\vbox{\hrule height .25pt width #4
\hbox to #4{\vrule width .25pt height #3\hfill\vrule width .25pt height #3}
\hrule height.25pt}\hss}
\vskip1sp
\setbox0=\hbox{\eightsc Figure #1.\hskip12pt\eightpoint\rm #2}
\ifdim\wd0>\hsize
\noindent\eightsc Figure #1.\hskip12pt\eightpoint\rm #2\else
\centerline{\eightsc Figure #1.\hskip12pt\eightpoint\rm #2}\fi
\vskip16pt
\endinsert}

\def\tabcaption#1#2{\vskip6pt
\setbox0=\hbox{\eightsc Table #1.\hskip12pt\eightpoint #2}
\ifdim\wd0>\hsize
\noindent\eightsc Table #1.\hskip12pt\eightpoint #2
\else
\centerline{\eightsc Table #1.\hskip12pt\eightpoint #2}
\fi
\vskip6pt}

\def\endinsert{\egroup\if@mid\dimen@\ht\z@\advance\dimen@\dp\z@ 
\advance\dimen@ 12\p@\advance\dimen@\pagetotal\ifdim\dimen@ >\pagegoal 
\@midfalse\p@gefalse\fi\fi\if@mid\smallskip\box\z@\bigbreak\else
\insert\topins{\penalty 100 \splittopskip\z@skip\splitmaxdepth\maxdimen
\floatingpenalty\z@\ifp@ge\dimen@\dp\z@\vbox to\vsize {\unvbox \z@ 
\kern -\dimen@ }\else\box\z@\nobreak\smallskip\fi}\fi\endgroup}

\def\pagecontents{
\ifvoid\topins \else\iftitle\else 
\unvbox \topins \fi\fi \dimen@ =\dp \@cclv \unvbox 
\@cclv 
\ifvoid\topins\else\iftitle\unvbox\topins\fi\fi
\ifvoid \footins \else \vskip \skip \footins \footnoterule 
\unvbox \footins \fi \ifr@ggedbottom \kern -\dimen@ \vfil \fi}


\newif\ifappend

\def\appendix#1#2{\def\applett{#1}\def\two{#2}%
\global\appendtrue
\global\theoremcount=0
\global\eqcount=0
\vskip18pt plus 18pt
\vbox{\parindent=0pt
\everypar={\hskip\parfillskip}
\def\\ {\vskip1sp}\elevenbold Appendix%
\ifx\applett\empty\gdef\applett{A}\ifx\two\empty\else.\fi%
\else\ #1.\fi\hskip6pt#2\vskip12pt}%
\global\sectiontrue%
\everypar={\global\sectionfalse\everypar={}}\nobreak\ignorespaces}

\newif\ifRefsUsed
\long\def\references{\global\RefsUsedtrue\vskip21pt%
\theinstitutions
\global\everypar={}\global\bibnum=0
\vskip20pt\goodbreak\bgroup%
\vbox{\centerline{\eightsc References}\vskip6pt}%
\ifdim\maxbibwidth>0pt
\leftskip=\maxbibwidth%
\parindent=-\maxbibwidth%
\else
\leftskip=18pt%
\parindent=-18pt%
\fi
\ninepoint%
\frenchspacing
\nobreak\ignorespaces\everypar={\amref}%
}

\def\endreferences{\vskip1sp\egroup\global\everypar={}%
\nobreak\vskip8pt\vbox{\thereceived\therevised}
}

\newcount\bibnum

\def\amref#1 {\global\advance\bibnum by1%
\immediate\write\auxfile{\string\expandafter\string\def\string\csname
\space #1croref\string\endcsname{[\the\bibnum]}}%
\leavevmode\hbox to18pt{\hbox to13.2pt{\hss[\the\bibnum]}\hfill}}

\def\bibline{\hbox to30pt{\hrulefill}\/\/}

\def\name#1{{\eightsc#1}}

\newdimen\maxbibwidth
\def\AuthorRefNames [#1] {%
\immediate\write\auxfile{\string\def\string\cite\string##1{[\string##1]}}

\def\amref{\spamref}
\setbox0=\hbox{[#1] }\global\maxbibwidth=\wd0\relax}

\def\spamref[#1] {\leavevmode\hbox to\maxbibwidth{\hss[#1]\hfill}}


\def\footnoterule{\kern-3pt\hrule width1in height.5pt\kern2.5pt}

\def\footnote#1#2{%
\plainfootnote{#1}{{\eightpoint\normalbaselineskip11pt
\normalbaselines#2}}}

\def\vfootnote#1{%
\insert \footins \bgroup \eightpoint\baselineskip11pt
\interlinepenalty \interfootnotelinepenalty
\splittopskip \ht \strutbox \splitmaxdepth \dp \strutbox \floatingpenalty 
\@MM \leftskip \z@skip \rightskip \z@skip \spaceskip \z@skip 
\xspaceskip \z@skip
{#1}$\,$\footstrut \futurelet \next \fo@t}


\newif\iffirstadded
\newif\ifadded

\def\addedlett{}

\def\alltheoremnums{%
\ifspecialnumon\global\specialnumonfalse
\ifadded\global\addedfalse
\iffirstadded\global\firstaddedfalse
\global\advance\theoremcount by1 \fi
\ifappend\applett\else\the\sectioncount\fi.\the\theoremcount\addedlett%
\xdef\theoremnum{\ifappend\applett\else\the\sectioncount\fi.%
\the\theoremcount\addedlett}%
\else$\rm\spnum$\def\theoremnum{{$\rm\spnum$}}\fi%
\else\global\firstaddedtrue
\global\advance\theoremcount by1 
\ifappend\applett\else\the\sectioncount\fi.\the\theoremcount%
\xdef\theoremnum{\ifappend\applett\else\the\sectioncount\fi.%
\the\theoremcount}\fi}

\def\allcorolnums{%
\ifspecialnumon\global\specialnumonfalse
\ifadded\global\addedfalse
\iffirstadded\global\firstaddedfalse
\global\advance\corolcount by1 \fi
\the\corolcount\addedlett%
\else$\rm\spnum$\def\corolnum{$\rm\spnum$}\fi%
\else\global\advance\corolcount by1 
\the\corolcount\fi}



\newif\ifthtitle

\let\saverparen)
\let\savelparen(
\def\rmparenl{{\rm(}}
\def\rmparenr{{\rm\/)}}
{
\catcode`(=13
\catcode`)=13
\gdef\makeparensRM{\catcode`(=13\catcode`)=13\let(=\rmparenl%
\let)=\rmparenr%
\everymath{\let(\savelparen%
\let)\saverparen}%
\everydisplay{\let(\savelparen%
\let)\saverparen\lookforbreak}}}

\medskipamount=8pt plus.1\baselineskip minus.05\baselineskip

\def\rmtext#1{\hbox{\rm#1}}

\def\proclaim#1{\vskip-\lastskip
\def\one{#1}\ifx\one\xtheorem\global\corolcount=0\fi
\ifsection\global\sectionfalse\vskip-6pt\fi
\medskip
{\elevensc#1}%
\ifx\one\xmaintheorem\global\corolcount=0
\gdef\theoremnum{Main Theorem}\else%
\ifx\one\xcorol\ \allcorolnums\else\ \alltheoremnums\fi\fi%
\ifthtitle\ \global\thtitlefalse{\rm(\thethtitle)}\fi.%
\hskip1em\bgroup\let\text\rmtext\makeparensRM\it\ignorespaces}

\def\nonumproclaim#1{\vskip-\lastskip
\def\one{#1}\ifx\one\xtheorem\global\corolcount=0\fi
\ifsection\global\sectionfalse\vskip-6pt\fi
\medskip
{\elevensc#1}.\ifx\one\xmaintheorem\global\corolcount=0
\gdef\theoremnum{Main Theorem}\fi\hskip.5pc%
\bgroup\it\makeparensRM\ignorespaces}

\def\endproclaim{\egroup\medskip}


\def\xproof{Proof}
\def\xremark{Remark}
\def\xcase{Case}
\def\xsubcase{Subcase}
\def\xconjecture{Conjecture}
\def\xstep{Step}
\def\xof{of}

\def\deconstruct#1 #2 #3 #4 #5 @{\def\one{#1}\def\two{#2}\def\three{#3}%
\def\four{#4}%
\ifx\two\empty #1\else%
\ifx\one\xproof%
\ifx\two\xof%
  \ifx\three\xcorol Proof of Corollary \rm#4\else%
     \ifx\three\xtheorem Proof of Theorem \rm#4\else\xone\fi%
  \fi\fi%
\else\xone\fi\fi.}

\def\pickup#1 {\def\this{#1}%
\ifx\this\xproof\global\let\go\demoproof
\global\let\enddemo\endproof\else
\ifx\this\xremark\global\let\go\demoremark\else
\ifx\this\xcase\global\let\go\demostep\else
\ifx\this\xsubcase\global\let\go\demostep\else
\ifx\this\xconjecture\global\let\go\demostep\else
\ifx\this\xstep\global\let\go\demostep\else
\global\let\go\demoproof\fi\fi\fi\fi\fi\fi}

\newif\ifnonum
\def\demo#1{\vskip-\lastskip
\ifsection\global\sectionfalse\vskip-6pt\fi
\def\one{#1 }\def\two{#1*}%
\setbox0=\hbox{\expandafter\pickup\one}\expandafter\go\two}

\def\numbereddemo#1{\vskip-\lastskip
\ifsection\global\sectionfalse\vskip-6pt\fi
\def\two{#1*}%
\expandafter\demoremark\two}

\def\demoproof#1*{\medskip\def\xone{#1}
{\ignorespaces\it\expandafter\deconstruct\xone {} {} {} {} {} @%
\unskip\hskip6pt}\rm\ignorespaces}

\def\demoremark#1*{\medskip
{\it\ignorespaces#1\/} \ifnonum\global\nonumtrue\else
 \alltheoremnums\unskip.\fi\hskip1pc\rm\ignorespaces}

\def\demostep#1 #2*{\vskip4pt
{\it\ignorespaces#1\/} #2.\hskip1pc\rm\ignorespaces}

\def\enddemo{\medskip}

\def\endproof{\ifmathqed\global\mathqedfalse\medskip\else
\parfillskip=0pt~~\hfill\qed\medskip
\fi\global\parfillskip0pt plus 1fil\relax
\gdef\enddemo{\medskip}}

\def\qed{\vbox{\hrule\hbox{\vrule height6pt\hskip6pt\vrule}\hrule}}


\def\proofbox{\parfillskip=0pt~~\hfill\qed\vskip1sp\parfillskip=
0pt plus 1fil\relax}







\def\stripbs#1#2*{\def\one{#2}}

\def\emptyspace{ }
\def\nextthing{}
\def\newline{***}
\def\eatone#1{ }

\def\lookatspace#1{\ifcat\noexpand#1\ \else%
\gdef\nextthing{}\xdef\next{#1}%
\ifx\next\emptyspace%
\let\nextthing\emptyspace\else\ifx\next\newline%
\gdef\nextthing{\eatone}\fi\fi\fi\egroup\nextthing#1}

{\catcode`\^^M=\active%
\gdef\spacer{\bgroup\catcode`\^^M=\active%
\let^^M=\newline\obeyspaces\lookatspace}}

\def\ref#1{\seeifdefined{#1}\expandafter\csname\one\endcsname\spacer}

\def\cite#1{\expandafter\ifx\csname#1croref\endcsname\relax[??]\else
\csname#1croref\endcsname\fi\spacer}


\def\seeifdefined#1{\expandafter\stripbs\string#1croref*%
\crorefdefining{#1}}

\newif\ifcromessage
\global\cromessagetrue

\def\crorefdefining#1{\ifdefined{\one}{}
{\ifcromessage\global\cromessagefalse%
\message{\spaces\spaces\spaces\spaces\spaces\spaces\spaces}%
\message{<Undefined reference.}%
\message{Please TeX file once more to have accurate cross-references.>}%
\message{\spaces\spaces\spaces\spaces\spaces\spaces\spaces}\fi[??]}}

\def\label#1#2*{\gdef\ctest{#2}%
\xdef\currlabel{\string#1croref}
\expandafter\seeifdefined{#1}%
\ifx\empty\ctest%
\xdef\labelnow{\write\auxfile{\noexpand\def\currlabel{\the\pageno}}}%
\else\xdef\labelnow{\write\auxfile{\noexpand\def\currlabel{#2}}}\fi%
\labelnow}

\def\ifdefined#1#2#3{\expandafter\ifx\csname#1\endcsname\relax%
#3\else#2\fi}




\def\articlecontents{
\vskip20pt\centerline{\bf Table of Contents}\everypar={}\vskip6pt
\bgroup \leftskip=3pc \parindent=-2pc 
\def\item##1{\vskip1sp\indent\hbox to2pc{##1.\hfill}}}

\def\endcontents{\vskip1sp\leftskip=0pt\egroup}

\def\journalcontents{\vfill\eject
\def\currannalsline{\hfill}
\global\titletrue
\vglue3.5pc
\centerline{\tensc\hskip12pt TABLE OF CONTENTS}\everypar={}\vskip30pt
\bgroup \leftskip=34pt \rightskip=-12pt \parindent=-22pt 
  \def\\ {\vskip1sp\noindent}
\def\pagenum##1{\unskip\parfillskip=0pt\dotfill##1\vskip1sp
\parfillskip=0pt plus 1fil\relax}
\def\name##1{{\tensc##1}}}


\institution{}
\onpages{0}{0}
\def\lastpage{???}
\def\thetitle{Title ???}
\def\theauthors{Authors ???}
\def\thereceived{}
\def\therevised{}

\gdef\split{\relaxnext@\ifinany@\let\next\insplit@\else
 \ifmmode\ifinner\def\next{\onlydmatherr@\split}\else
 \let\next\outsplit@\fi\else
 \def\next{\onlydmatherr@\split}\fi\fi\let\eqnu\xspliteqnu\next}

\gdef\align{\relaxnext@\ifingather@\let\next\galign@\else
 \ifmmode\ifinner\def\next{\onlydmatherr@\align}\else
 \let\next\align@\fi\else
 \def\next{\onlydmatherr@\align}\fi\fi\let\eqnu\xspliteqnu\next}

\def\spliteqnu{{\tenrm\sectandeqnum}\relax}

\def\xspliteqnu{\tag\spliteqnu}

\catcode`@=12

\document


\annalsline{June}{1997}
\startingpage{1}     


%
%
%
%
%

\def\for{\  \hbox{ for } \ }
\def\if{ \ \hbox{ if } \ }

\def\where{\  \hbox{ where } \ }
\def\and{\  \hbox{ and } \ }

\def\equal{\buildrel  def \over =}

\def\la{\lambda}
\def\La{\Lambda}
\def\om{\omega}

\def\th{\theta}
\def\al{\alpha}

\def\de{\delta}
\def\De{\Delta}
\def\ka{\kappa}
\def\si{\sigma}

\def\Ga{\Gamma}
\def\ze{\zeta}

\def\pa{\partial}

\def\vph{\varphi}

\def\vep{\varepsilon}

\def\tal{\tilde{\alpha}}

\def\tw{\tilde w}

\def\tz{\tilde z}
\def\tb{\tilde b}

\def\tG{\tilde G}

\def\hH{\hat{H}}

\def\hY{\hat{Y}}

\def\hT{\hat{T}}

\def\hw{\hat{w}}
\def\hu{\hat{u}}

\def\hv{\hat{v}}

\def\C{\bold{C}}

\def\R{\bold{R}}

\def\Z{\bold{Z}}

\def\one{\bold{1}}

\def\0{\bold{0}}

\def\C{\hbox{\bf C}}


\def\f{\Cal{F}}

\def\l{\Cal{L}}

\def\g{\Cal{G}}

\def\o{\Cal{O}}

\font\germ=eufb10 
\def\goth#1{\hbox{\germ #1}}

\def\TT{\goth{T}}
\def\HH{\goth{H}}

\def\AA{\goth{A}}

\font\smm=msbm10 at 12pt 
\def\symbol#1{\hbox{\smm #1}}
\def\lsmash{{\symbol n}}



\title
{Inverse Harish-Chandra transform \\
and difference operators}
 
\shorttitle{ Harish-Chandra transform }

\acknowledgements{
Partially supported by NSF grant DMS--9622829}

\author{ Ivan Cherednik}

\institutions{
Math. Dept, University of North Carolina at Chapel Hill,   
 N.C. 27599-3250
\\ Internet: chered\@math.unc.edu
}

%
%
%
%
%
\vfil

In the paper we calculate the images of
the operators of multiplication by Laurent polynomials with respect
to the Harish-Chandra transform and its non-symmetric generalization
due to Opdam. It readily leads to  a new simple proof
of the Harish-Chandra inversion theorem
in the zonal case 
(see [HC,He1]) and the corresponding theorem from [O1].
We assume that $k> 0$ 
and restrict ourselves to compactly
supported functions, borrowing the growth estimates
from [O1]. 

The Harish-Chandra transform is the integration of 
symmetric functions
on the maximal real split torus $ A$ of a 
semi-simple Lie group $G$
multiplied by  the spherical zonal function $\phi(X,\la)$, where
$X\in A,\ \la \in (\hbox{Lie}A\otimes_{\R}\C)^*.$
The measure is the restriction of the invariant measure on $G$
to the space of double cosets $K\backslash G/K \subset A/W$ for  the maximal
compact subgroup $K\subset G$ and  the restricted  Weyl group $W$.
The function $\phi$ is a symmetric ($W$-invariant) eigenfunction
of the radial parts of the $G$-invariant differential operators on $G/K$; 
$\la$ determines the set of eigenvalues. 
The parameter $k$ is given by the root multiplicities 
($k=1$ in the group case). There is a generalization to arbitrary $k$ 
due to Calogero, Sutherland, Koornwinder, Moser, Olshanetsky, Perelomov,
Heckman, and Opdam. See [HO1,H1,O2] for a systematic theory.  
In the non-symmetric variant due to Opdam [O1], the operators from
[C5] replace the radial parts of $G$-invariant operators  and their
$k$-generalizations. The problem is to define the inverse transforms
for various classes of functions.

In the papers [C1-C4], a difference counterpart of the Harish-Chandra
transform was suggested, which is also a deformation of the
Fourier transform in the $p$-adic theory of spherical functions. 
Its  kernel (a $q$-generalization of $\phi$)
is defined as  an eigenfunction
of the $q$-difference ``radial parts'' (Macdonald's operators and their
generalizations). 
There are  applications 
in combinatorics (the Macdonald polynomials), 
the representation theory (say, quantum groups at roots of unity),
and the mathematical physics 
(the Knizhnik-Zamolodchikov equations and more). 
The difference Fourier transform  is  self-dual, i.e. the kernel
is $x\leftrightarrow \la$ symmetric for $X=q^x$.
This holds in the differential theory only 
for the so-called rational degeneration with the tangent space 
$T_e(G/K)$ instead of $G/K$ (see [He1,D,J]) and for special $k=0,1$.
In the latter cases, the differential and
difference transforms coincide up to a normalization.

The differential case is recovered when  $q$ goes to $1$.
At the moment the analytic methods are not mature enough to 
manage the limiting procedure in detail. So
we develop the corresponding
technique  without any reference to the difference Fourier transform, 
which is not introduced and discussed in the paper, switching from 
the double affine Hecke algebra [C9] to  
the degenerate one from [C7]. The latter generalizes Lusztig's graded
affine Hecke algebra [L] (the $GL_n$-case is due to Drinfeld).
Mainly we need the intertwiners and creation operators
from [C2] (see [C9] and [KS] for $GL_n$). 
It is worth mentioning
that the $p$-adic theory  corresponds to the limit $q\to 0$,
which does not violate the standard difference
convergence condition $|q|<1$ and hardly creates any analytic problems.

The key result is that 
the Opdam transforms  of the operators of 
multiplication by the coordinates
coincide with the operators from  
[C6] (see (\ref\Desb) below).
Respectively, the Harish-Chandra transforms of 
the multiplications by $W$-invariant Laurent polynomials
are the {\it difference} operators from (\ref\Laop).
It is not very surprising
that the transforms of these important operators
haven't been found before.
The calculation, not difficult by itself,
involves  the following ingredients new in the
harmonic analysis on symmetric spaces:

\vfill
\item{ a)} The Dunkl-type operators [C5],
their relation  to degenerate (gra\-ded) affine 
Hec\-ke algebras [C8], and the Opdam transform [O1].

\item{ b)} The Macdonald operators,
difference counterparts of  Dunkl operators,
and double affine Hecke algebras  
[M1,M2,C6,C9].
\vfill

We begin the paper with some basic facts about double affine
Hecke algebras and their polynomial representations which are necessary
to produce dif\-fe\-rence Dunkl-type operators. Then we 
discuss the two degeneration procedures  leading respectively
to the operators from [C5](differential) and from [C6](difference).
Both are governed by {\it the same} degenerate double affine
Hecke algebra (introduced in [C7]). It eventually
results in the inversion theorem.

Hopefully  this approach can be extended  to any $k$
and, moreover, may help to (re)establish the
analytic properties of the direct and inverse
transforms. 
I believe that this paper will stimulate a {\it systematic} renewal
of the Harish-Chandra
theory and related representation theory on the basis of the
difference-operator methods.

The paper was written during my stay at RIMS
(Kyoto University).
I am grateful 
for the kind invitation and  hospitality.
I thank M. Kashiwara for valuable remarks and
 M. Duflo, J. Faraut, D. Kazhdan,  T. Miwa, E. Mukhin,
and E. Opdam for useful discussions.

%
%
%
%
\section {Affine Weyl groups}

Let $R=\{\al\}   \subset \R^n$ be a root system of type $A,B,...,F,G$
with respect to a euclidean form $(z,z')$ on $\R^n \ni z,z'$,
normalized by the standard condition
that $(\al,\al)=2$ for long $\al$.
Let us  fix the set $R_{+}$ of positive  roots,
the corresponding simple 
roots $\al_1,...,\al_n$, and  their dual counterparts 
$a_1 ,..., a_n,  a_i =\al_i^\vee, \where \al^\vee =2\al/(\al,\al)$.  
The dual fundamental weights
$b_1,...,b_n$  are determined by the relations  $ (b_i,\al_j)= 
\de_i^j $ for the 
Kronecker delta. We will also use the dual root system
$R^\vee =\{\al^\vee, \al\in R\}, R^\vee_+$, and 
$$
\eqalignno{
& A=\oplus^n_{i=1}\Z a_i \subset B=\oplus^n_{i=1}\Z b_i,\ 
B_+=\oplus^n_{i=1}\Z_+ b_i\for \Z_+=\{m\ge 0\}. 
}
$$
In the standard notations, $A= Q^\vee,\ 
B = P^\vee $ -- see [B].  Later on,
$$
\eqalign{
&\nu_{\al}\  =\ (\al,\al),\  \nu_i\ =\ \nu_{\al_i}, \ 
\nu_R\ = \{\nu_{\al}, \al\in R\}\subset \{2,1,2/3\}.
}
$$

The vectors $\ \tal=[\al,k] \in 
\R^n\times \R = \R^{n+1}$ 
for $\al \in R, k \in \Z $ 
form the {\it affine root system} 
$R^a \supset R$ ( $z\in \R^n$ are identified with $ [z,0]$).
We add  $\al_0 \equal [-\th,1]$ to the  simple roots 
for the {\it maximal root} $\th \in R$.
The corresponding set $R^a_+$ of positive roots coincides
with $R_+\cup \{[\al,k],\  \al\in R, \  k > 0\}$. Let $a_0=\al_0$. 

We will use the affine Dynkin diagram $\Ga^a$  with
$\{\al_j,0 \le j \le n\}$ as the vertices  
($m_{ij}=2,3,4,6$\  if $\al_i\and\al_j$ are joined by 0,1,2,3 laces
respectively).
The set of
the indices of the images of $\al_0$ by all 
the automorphisms of $\Ga^a$ will be denoted by $O$ ($O=\{0\} 
\for E_8,F_4,G_2$). Let $O^*=\{r\in O, r\neq 0\}$.
The elements $b_r$ for $r\in O^*$ are the so-called minuscule
coweights ($(b_r,\al)\le 1$ for
$\al \in R_+$).

Given $\tal=[\al,k]\in R^a,  \ b \in B$, let  
$$
\eqalignno{
&s_{\tal}(\tz)\ =\  \tz-(z,\al^\vee)\tal,\ 
\ b'(\tz)\ =\ [z,\ze-(z,b)]
&\eqnu
\label\saction\eqnum*
}
$$
for $\tz=[z,\ze] \in \R^{n+1}$.

The {\it affine Weyl group} $W^a$ is generated by all $s_{\tal}$
(simple reflections $s_j=s_{\al_j} \for 0 \le j \le n$
are enough).
It is
the semi-direct product $W\lsmash A$, where the non-affine Weyl  
group $W$ is the span of $s_\al,
\al \in R_+$.
Here and further we identify $b\in B$ with the corresponding
translations. For instance,
$$
\eqalignno{
& \al^\vee =\ s_{\al}s_{[\al,1]}=\ s_{[-\al,1]}s_{\al} \for
\al\in R.
&\eqnu
\label\baction\eqnum*
}
$$

The {\it extended Weyl group} $ W^b$ generated by $W\and B$
 is isomorphic to $W\lsmash B$:
$$
\eqalignno{
&(wb)([z,\ze])\ =\ [w(z),\ze-(z,b)] \for w\in W, b\in B.
&\eqnu
}
$$
For $b_+\in B_+$, let
$$
\eqalignno{
&\om_{b_+} = w_0w^+_0  \in  W,\ \pi_{b_+} =
b_+(\om_{b_+})^{-1}
\ \in \ W^b, \ \om_i=\om_{b_i},\pi_i=\pi_{b_i},
&\eqnu
\label\wo\eqnum*
}
$$
where $w_0$ (respectively, $w^+_0$) is the longest element in $W$
(respectively, in $ W_{b_+}$ generated by $s_i$ preserving $b_+$) 
relative to the 
set of generators $\{s_i\}$ for $i >0$.

The elements $\pi_r\equal\pi_{b_r}, r \in O^*,$ and $\pi_0=\hbox{id}$
 leave $\Ga^a$ invariant 
and form a group denoted by $\Pi$, 
 which is isomorphic to $B/A$ by the natural 
projection $b_r \mapsto \pi_r$. As to $\{\om_r\}$,
they preserve the set $\{-\th,\al_i, i>0\}$.
The relations $\pi_r(\al_0)= \al_r= (\om_r)^{-1}(-\th) 
$ distinguish the
indices $r \in O^*$. Moreover (see e.g. [C6]):
$$
\eqalignno{
& W^b  = \Pi \lsmash W^a, \where
  \pi_rs_i\pi_r^{-1}  =  s_j \if \pi_r(\al_i)=\al_j,\  0\le i,j\le n.
&\eqnu
}
$$
The {\it length} $l(\hw)$  of 
$\hw = \pi_r\tw \in W^b$ is the length of 
a reduced  
decomposition $\tw\in W^a$ with respect to
$\{s_j, 0\le j\le n\}$:
$$
\eqalign{
l(\hw)=|\ell(\hw)|,\ \ell(\hw)\ & =\  R^a_+ \cap \hw^{-1}(-R^a_+)
 \cr
& =\ \{\tal\in R^a_+, \ l( \hw s_{\tal}) < l(\hw) \}.
}
 \eqnu
\label\lambda\eqnum*
$$
For instance, given $b_+\in B_+$,
$$
\eqalignno{
&\ell(b_+)\  =\ \{\tal= [\al,k],\ \ \al\in R_+\and    
( b_+, \al )>k\ge 0 \},\cr
& l(b_+)\ =\ 2(b_+,\rho),  
\where \rho=(1/2)\sum_{\al\in R_+}\al. 
&\eqnu 
\label\lambi\eqnum* 
}
$$
We will also use  the dominant {\it affine Weyl chamber} 
$$
\eqalignno{
&C^a_+\ =\ \{ z\in \R^n \hbox{ \ such\ that\  } (z,\al_i)>0 \for i>0,\
(z,\th)<1\}.
&\eqnu
\label\afchamber\eqnum*
}
$$
%
%
%
%
%
\section{ Double affine Hecke algebras} 
We put    
$m=2 \for D_{2l} \and C_{2l+1},\ m=1 \for C_{2l}, B_{l}$,
otherwise $m=|\Pi|$. We consider $ q,\ \{ t(\nu) , \nu \in \nu_R \},\ $
 $X_1,\ldots,X_n$
 as  independent variables, setting
$$
\eqalignno{
&   t_{\tal} = t(\nu(\al)),\ t_j = t_{\al_j},
\where \tal=[\al,k] \in R^a,\ 0\le j\le n, \cr 
& X_{\tb}\ =\ \prod_{i=1}^nX_i^{k_i} q^{ k} 
\if \tb=[b,k] \for b=\sum_{i=1}^nk_i b_i\in B.
&\eqnu  
\label\Xde\eqnum*
}
$$

Later on $ k\in (1/ 2m)\Z$,  
$ \C_{ q,t}$
 is the field of rational 
functions in terms of $ q^{1/2m},
\{t(\nu)^{1/2} \}$,
$\C_{q,t}[X] = \C_{q,t}[X_b]$  means the algebra of
polynomials in terms of $X_i^{\pm 1}$ 
with the coefficients from $ \C_{ q,t}$.
We set  $\al_{r^*} \equal \pi_r^{-1}(\al_0)$, so
$\om_r\om_{r^*}=1=\pi_r\pi_{r^*}$.

\proclaim{Definition [C9]}
 The  double  affine Hecke algebra $\HH\ $
is generated over the field $ \C_{ q,t}$ by 
the elements $\{ T_j,\ 0\le j\le n\}$, 
pairwise commutative $\{X_b, \ b\in B\}$ satisfying (\ref\Xde),
 and the group $\Pi$ where the following relations are imposed:

(o)\ \  $ (T_j-t_j^{1/2})(T_j+t_j^{-1/2})\ =\ 0,\ 0\ \le\ j\ \le\ n$;

(i)\ \ \ $ T_iT_jT_i...\ =\ T_jT_iT_j...,\ m_{ij}$ factors on each side;

(ii)\ \   $ \pi_rT_i\pi_r^{-1}\ =\ T_j \if \pi_r(\al_i)=\al_j$; 

(iii)\  $T_iX_b T_i\ =\ X_b X_{a_i}^{-1} \if (b,\al_i)=1,\
1 \le i\le  n$;

(iv)\  $T_0X_b T_0\ =\ X_{s_0(b)}\ =\ X_b X_{\th} q^{-1}
\if (b,\th)=-1$;

(v)\ \ $T_iX_b\ =\ X_b T_i \if (b,\al_i)=0 (i>0),\ 
(b,\th)=0 (i=0)$;

(vi)\ $\pi_rX_b \pi_r^{-1}\ =\ X_{\pi_r(b)}\ =\ X_{\om^{-1}_r(b)}
 q^{(b_{r^*},b)},\  r\in O^*$.
\label\double\theoremnum*
\endproclaim

Given $\tw \in W^a, r\in O,\ $ the product
$$
\eqalignno{
&T_{\pi_r\tw}\equal \pi_r\prod_{k=1}^l T_{i_k},\where 
\tw=\prod_{k=1}^l s_{i_k},
l=l(\tw),
&\eqnu
\label\Tw\eqnum*
}
$$
does not depend on the choice of the reduced decomposition
(because $\{T\}$ satisfy the same ``braid'' relations as $\{s\}$ do).
Moreover,
$$
\eqalignno{
&T_{\hv}T_{\hw}\ =\ T_{\hv\hw}\  \hbox{ whenever}\ 
 l(\hv\hw)=l(\hv)+l(\hw) \for
\hv,\hw \in W^b.
&\eqnu}
$$
\label\TT\eqnum*
  In particular, we arrive at the pairwise 
commutative elements 
$$
\eqalignno{
& Y_{b}\ =\  \prod_{i=1}^nY_i^{k_i} \if  
b=\sum_{i=1}^nk_ib_i\in B,\where  
 Y_i\equal T_{b_i},
&\eqnu
\label\Yb\eqnum*
}
$$
satisfying the relations
$$
\eqalign{
&T^{-1}_iY_b T^{-1}_i\ =\ Y_b Y_{a_i}^{-1} \if (b,\al_i)=1,
\cr
& T_iY_b\ =\ Y_b T_i \if (b,\al_i)=0, \ 1 \le i\le  n.}
\eqnu
$$

\comment 
The map
$$
\eqalign{
 \vph: &X_i \mapsto Y_i^{-1},\ \  Y_i \mapsto X_i^{-1},\  \ T_i \mapsto T_i, \cr
&t(\nu) \mapsto t(\nu),\
 q\mapsto  q,\ \nu\in \nu_R,\ 1\le i\le n.
}
\eqnu
\label\vph\eqnum*
$$
can be extended to an anti-involution 
($\vph(AB)=\vph(B)\vph(A)$)
 of \HH\ .
The following maps can be extended to involutions of
\HH\ (see [C1,C3]):
$$
\eqalignno{
  \vep:\ &X_i \mapsto Y_i,\ \  Y_i \mapsto X_i,\  \ T_i \mapsto T_i^{-1},
&\eqnu
\label\vep\eqnum*
 \cr
&t(\nu) \mapsto t(\nu)^{-1},\
 q\mapsto  q^{-1},\cr
 \tau: \ &X_b \mapsto X_b,\ \ Y_r \mapsto X_rY_r q^{-(b_r,b_r)/2},\ 
Y_\th \mapsto X_0^{-1}T_0^{-2}Y_\th, &\eqnu
\cr
&T_i\mapsto T_i,\ t(\nu) \mapsto t(\nu),\
 q\mapsto  q, \cr
&1\le i\le n,\ r\in O^*,\ X_0=qX_\th^{-1}.
\label\tau\eqnum*
}
$$
Let us give some explicit formulas:
$$
\eqalign{
&\vep(T_0)\ =\ X_\th T_0^{-1}Y_\th\ =\ X_\th T_{s_\th},\ 
\vep(\pi_r)\ =\ X_rT_{\om_r^{-1}},\cr
&\tau(T_0)\ =\ X_0^{-1}T_0^{-1},\ \tau(\pi_r)\ =
\ q^{-(b_r,b_r)/2}X_r\pi_r= q^{(b_r,b_r)/2}\pi_rX_{r^*}^{-1},\cr
&\pi_r X_{r^*}\pi_r^{-1}\ =\ q^{(b_r,b_r)}X_{r}^{-1},\ 
X_{r^*}T_{\om_r}X_r\ =\ T_{\om_{r^*}}^{-1}.
}
\eqnu
\label\vphto\eqnum*
$$
The  involution $\eta=\tau^{-1}\vep\tau$
$$
\eqalignno{
  \eta: \ &X_r\mapsto  q^{(b_r,b_r)/2}X_r^{-1}Y_r\ =\  q^{-(b_r,b_r)/2}
\pi_rX_{r^*}T_{\om_r},
&\eqnu
\cr
 &Y_r \mapsto q^{(b_r,b_r)}X_r^{-1}Y_rX_r\  
 =\ \pi_rT_{\om_{r^*}}^{-1},\cr 
&Y_\th \mapsto T_0^{-1}T_{s_\th}^{-1},\ 
T_j\mapsto T_j^{-1} (0\le j\le n),\cr
&\pi_r\mapsto\pi_r (r\in O^*),\ t(\nu) \mapsto t(\nu)^{-1},\
 q\mapsto  q^{-1}.
\label\eta\eqnum*
}
$$
We note that $\vep$ and $\eta$ commute with 
the main anti-involution $^*$ from [C2]:
\endcomment

The following {\it basic anti-involution} $^*$ [C11]
plays the key role for the so-called
polynomial (basic) representation
of \HH\  and difference counterparts
of the Dunkl operators:
$$
\eqalign{
  & X_i^*\ =\  X_i^{-1},\ \  Y_i^*\ =\  Y_i^{-1},\  \
 T_i^* \ =\  T_i^{-1},\ \pi_r^*\ =\ \pi_r^{-1},\cr
&t(\nu)^*= t(\nu)^{-1},\
 q^*= q^{-1},\ 0\le i\le n,\ 
(AB)^*=B^*A^*,\ A,B\in \HH\ .
}
\eqnu
\label\star\eqnum*
$$

The {\it $Y$-intertwiners} (see  [C2])
are introduced as follows:
$$
\eqalignno{
&\Phi_j\ =\ 
T_j + (t_j^{1/2}-t_j^{-1/2})(Y_{a_j}^{-1}-1)^{-1},\ 1\le j\le n,
\cr
&\Phi_0\ =\ 
X_\th T_{s_\th} - (t_0^{1/2}-t_0^{-1/2})(q^{-1}Y_{\th}^{-1}-1)^{-1},
\cr
&P_r\ =\ 
X_rT_{\om_r^{-1}},\ r\in O^*. 
\label\Phi\eqnum*
}
$$
\comment
& G_j=\Phi_j (\phi_j)^{-1},\
\tG_j=(\phi_j)^{-1}\Phi_j,\ 
 &\eqnu
\cr 
&\phi_j\ =\  t_j^{1/2} + 
(t_j^{1/2} -t_j^{-1/2})(X_{a_j}-1)^{-1},
\label\Phi\eqnum*
}
$$
for $0\le j\le n$.
\endcomment

 They belong to the 
extension \HH\ by the field $\C_{q,t}(Y)$
of rational functions in $\{Y\}$. The elements
$\{\Phi_j,\ P_r\}$ satisfy the relations 
for $\{T_j, \pi_r\}$.  Hence the elements 
$$
\Phi_{\hw}\ =\ P_r \Phi_{j_l}\cdots \Phi_{j_1},
\where \hw=\pi_r s_{j_l}\cdots s_{j_1}\in W^b,
\eqnu
\label\Phiprod\eqnum*
$$
are well-defined 
for reduced decompositions of $\hw$,
and $\Phi_{\hw}\Phi_{\hu}=\Phi_{\hw\hu}$ whenever
$l(\hw\hu)=l(\hw)+l(\hu)$.

The following 
property of $\{\Phi\}$  fixes them uniquely up to left
or right multiplications by functions of $Y$:
$$
\eqalignno{
&\Phi_{\hw} Y_b\ =\ Y_{\hw(b)}\Phi_{\hw},\ \hw\in W^b,
&\eqnu 
\label\Phix\eqnum*
}
$$
where $Y_{[b,k]}\equal Y_b q^{-k}$.

The  $\{\Phi_{\hw}\}$ are exactly the images of the 
$X$-intertwiners (see [C2], (2.12)) with respect to
the involution $\vep$ from (2.6) ibidem. 

Note that 
$$
\eqalignno{
& \Phi_{\hw}^* = \Phi_{\hw^{-1}},  \ \hw\in W^b.
&\eqnu 
\label\Phistar\eqnum*
}
$$
Use the quadratic relation for $T_j$ to check it for 
$\{\Phi_j,0\le j\le n\}$.

In the case of $GL_n$, the element $P_1$ (it is of infinite
order) was used by Knop and Sahi (see [Kn]) to construct the
non-symmetric Macdonald polynomials and eventually confirm
the Macdonald integrality conjecture. They introduced it
independently of [C9] and checked the ``intertwining''
property directly. Generally speaking, a direct verification
of (\ref\Phix) without the $X$-intertwiners and the 
involution $\vep$  
is possible but much more difficult.

%
%
%
%
\section { Polynomial representation} 
We will identify $X_{\tb}$ with the corresponding multiplication
operators: 
$$
\eqalignno{
& X_{\tb} (p(X))\ =\ X_{\tb} p(X),\    p(X) \in
\C_{q,t} [X],
&\eqnu}
$$
\label\X\eqnum*
and  use the action (\ref\saction):
$$
\eqalignno{
&\hw(X_{b})\ =\ X_{\hw(b)},\hbox{\ for\ instance,\ }
\cr
&s_0(X_b) = X_b X_{\th}^{-(b,\th)} q^{(b,\th)},\ 
a(X_b)=q^{-(a,b)}X_b \for a,b\in B.
&\eqnu
\label\afxact\eqnum*
}
$$

The {\it Demazure-Lusztig operators}
$$
\eqalignno{
&\hT_j\  = \  t_j ^{1/2} s_j\ +\ 
(t_j^{1/2}-t_j^{-1/2})(X_{a_j}-1)^{-1}(s_j-1),
\ 0\le j\le n.
&\eqnu
\label\Demaz\eqnum*
}
$$
are well-defined on $\C_{ q,t}[X]$.

We note that only $\hT_0$ depends on $ q$: 
$$
\eqalign{
&\hT_0\  =\  t_0^{1/2}s_0\ +\ (t_0^{1/2}-t_0^{-1/2})
( q X_{\th}^{-1} -1)^{-1}(s_0-1).
}
\eqnu
$$

\proclaim{Theorem [C9]}
 The map $ T_j\mapsto \hT_j,\ X_b \mapsto X_b$,
$\pi_r\mapsto \pi_r$  can be extended to a
faithful representation $H\mapsto\hH$ of
\HH\ (fixing $q,t$). 
The operators $\l_p=p(\hY_1,\ldots,\hY_n)$ are
$W$-invariant for $W$-invariant polynomials $p$ 
and preserve $\C_{q,t}[X]^W$. 
\endproclaim
\proofbox
\label\faith\theoremnum*

The exact formulas for the restrictions $L_p$ of
$\l_p$ to $W$-invariant functions are known only for rather simple $p$.

\proclaim{Proposition [C11]}
Given $r\in O^*$, let $m_r = \sum _{w \in W/W_r} X_{w(-b_r)}$, where $W_r$ is
the stabilizer of $b_r$ in $W$. Then $\ell_r=\ell(b_r)\subset R_+$ and
$$
\eqalign{
& L_r=L_{m_r} =\ \sum_{w\in W/W_r}
 \prod _{\al \in \ell_r}
{t^{1/2}_\al X_{w(\al^\vee)}- t^{-1/2}_{\al}\over
X_{w(\al^\vee)}- 1} \, w(-b),
}
\eqnu
\label\macoper\eqnum*
$$
where $w(-b)=-w(b)\in B$ acts as in (\ref\afxact). 
\label\MACOPER\theoremnum*
\endproclaim
\proofbox

The operators $L_r$ were introduced by Macdonald [M1,M2].
For $GL_n$, they were found independently by Ruijsenaars [R].
Later we will use these formulas in the ``rational'' limit.

The coefficient of $X^0=1$ ({\it the constant term})
of a polynomilal
$f\in  \C_{q,t}[X]$
will be denoted by $\langle  f \rangle$. Let
$$
\eqalign{
&\mu\ =\ \prod_{a \in R_+^\vee}
\prod_{i=0}^\infty {(1-X_aq_a^{i}) (1-X_a^{-1}q_a^{i+1})
\over
(1-X_a t_aq_a^{i}) (1-X_a^{-1}t_a^{}q_a^{i+1})},
}
\eqnu
\label\mu\eqnum*
$$
where $q_a=q^{(a,a)/2}=q^{2/(\al,\al)} \for a=\al^\vee$.
 
The coefficients of $\mu_0\equal \mu/\langle \mu \rangle_0$
are from $\C(q,t)$, and
$\mu_0^*\ =\ \mu_0$  with respect to the involution 
$$
 X_b^*\ =\  X_{-b},\ t^*\ =\ t^{-1},\ q^*\ =\ q^{-1}.
$$
Hence
$$
\eqalignno{
&\{f,g\}_\mu\equal \langle  f\, {g}^*\mu_0\rangle\ =\ 
(\{g,f\}_\mu)^* \for
f,g \in \C(q,t)[X].
&\eqnu
\label\innerpro\eqnum*  
}
$$
Here and further see [C3,C11].

\proclaim {Proposition}
 For any $H\in \HH\ $ and the anti-involution $^*$
from (\ref\star), $\{ \hH(f), g\}_\mu  = $ 
$ \{ f, \hH^*(g)\}_\mu$ for   $f,g\in \C_{q,t}[X]$.
The operators $\{X_b,Y_b, T_j,
\pi_r\}$ and $q,t$ are unitary with respect to the form $\{\ ,\ \}_\mu.$
\endproclaim
\label\YONE\theoremnum*
\proofbox

Actually this proposition is entirely algebraic.
It  holds for
other inner products. 
For instance, the inner product in terms of
Jackson integrals  considered in [C2] can be taken, which
is expected to have applications  to negative $k$.

%
%
%
%
\section { Degenerations} 

Let us fix 
$\ka=\{ \ka(\nu) \in \C ,\ \nu \in \nu_R\},\ \ka_\al=\ka((\al,\al)),\ 
\ka_j=\ka_{\al_j}$  and
introduce
$ \rho_\ka= (1/2)\sum _{\al\in R_+}\ka_\al \al.$
One has: $(\rho_\ka,\al_j^\vee)=\ka_j$ for $1\le j\le n$.

\proclaim {Definition}
The   degenerate (graded) double affine Hecke algebra $\HH' $ is 
generated by
the group algebra $\C [W^b]$ and the pairwise commutative
$$
\eqalignno{
&y_{\tb}\equal
\sum ^n_{i=1}(b,\al_i)y_i +u \for 
\tb=[b,u]\in \R^n\times\R,
&\eqnu
}
$$ 
satisfying  the following relations:
$$
\eqalignno{
&s_j y_{\tb}-y_{s_j(\tb)}s_j\ =\ - \ka_j(b,\al_j ),\ (b,\al_0)=-(b,\th),
\cr
&\pi_r y_{\tb}\ =\ y_{\pi_r(\tb)}\pi_r \for 
  0\le j\le n, \  r\in O.
&\eqnu
\label\suka\eqnum*
}
$$
\label\DEGHE\theoremnum*
\endproclaim

Without $s_0$ and $\pi_r$, we arrive at the defining relations
of the graded affine Hecke algebra from  [L] (see also [C8]).
It is a natural degeneration  
of the double affine
Hecke algebra when $q\to 1, t\to 1$.

We will also  use the parameters
$k_\al\equal(\al,\al)\ka_\al/2,\ k_i=k_{\al_i},$
and the derivatives of $\C[X]$: 
$$
 \partial_{a}(X_{b})\ =\ (a ,b)X_{b},\ a,b\in B.  
$$
Note that 
$w(\partial _{b})=
\partial _{w(b)}, \  w\in W$.

{\it Differential degeneration.}
 Following [C2,C10], let
$$
\eqalignno{
&q=\exp(h),\ t_j=q_j^{k_j}=q^{\ka_j}, 
\ Y_b=1-h y_b,
&\eqnu
\label\onelim\eqnum*
}
$$ 
and $h$ tend to zero. The terms of order $h^2$ are ignored
(without touching $X_b,\ka_j$).
We will readily arrive at the relations
(\ref\suka) for $y_b$.

Setting $X_b=e^{x_b}$, the limit of $\mu$ from (\ref\mu)  is 
$$
\eqalignno{
&\tau\equal \prod_{\al\in R_+}|2\sinh(x_{\al^\vee}/2)|^{2k_\al}.
&\eqnu
\label\tauform\eqnum*
}
$$

\proclaim{ Proposition}
a) The following operators acting on Laurent polynomials $f\in \C[X]$ 
$$
\eqalignno{
D_{b} \equal
&\pa_b +
\sum_{\al\in R_+} { \ka_{\al}(b,\al)\over
(1-X_{\al^\vee}^{-1}) }
\bigl( 1-s_{\al} \bigr)- (\rho_{\ka},b)
&\eqnu
\label\dunk\eqnum*
}
$$
are pairwise commutative, and $y_{[b,u]}=D_b+u$
satisfy (\ref\suka) for 
the following action of the group $W^b$:
$$ w^x(f)=w(f)\for w\in W,\ b^x(f)=X_b f \for b\in B.$$
For instance, $s_0^x(f)=X_\th s_\th(f),\ \pi_r^x(f)=X_r \om_r^{-1}(f).$

b) The operators $D_b$ are formally self-adjoint with respect to the inner
product 
$$
\eqalignno{
&\{f,g\}_\tau\equal \int f(x)g(-x)\tau dx,
&\eqnu
\label\tauprod\eqnum*
}
$$
i.e. $\tau^{-1}D_b^+\tau=D_b$ for the anti-involution $^+$
sending $\hw^x\mapsto (\hw^x)^{-1},\ \partial_b\mapsto -\partial_b,$ 
where $\hw\in W^b$.
\label\DUNDEG\theoremnum*
\endproclaim
{\it Proof.}
First, $y_b=\lim_{h\to 0}(1-\hY_b)/h$
is precisely  (\ref\dunk). Degenerate the product formulas
for $Y_i=T_{b_i}$ in the polynomial representation (cf. [C6,C11],
and [C10]). 
As to the term $-(\rho_\ka,b)$, 
use (\ref\lambi).
The condition
$q^*=q^{-1}$ formally leads to $h^*=-h$, and $k^*=k$. 
Similarly, $x_b^*$ must be equal to $-x_b$, and all these
give that $D_b^*=D_b$ for the inner product
(\ref\tauprod). Use Proposition \ref\YONE.
Actually the claim can be checked directly
without the limiting procedure. However we want to demonstrate that
the pairing (\ref\tauprod) is an exact counterpart of the
difference one. We mention that Opdam uses  somewhat different
pairing and involution in [O1].
\proofbox

Degenerating $\{ \Phi\}$, we get the intertwiners of 
$\HH'\ $:
$$
\eqalignno{
&\Phi'_i\ =\ 
s_i + {k_i \over y_{\al_i} },\ 
\Phi'_0\ =\ 
X_\th s_\th + {k_0 \over 1-y_{\th} },\cr
&P_r'\ =\ X_r\om_r^{-1}, \for  1\le i\le n,\ r\in O^*.
&\eqnu
\label\Phiprime\eqnum*
}
$$
The operator $P_1'$ in the case of $GL_n$ (it is of 
infinite order) plays the key role in [KS].
The formulas for $\Phi'_i$ when $1\le i\le n$ are well-known
in the theory of degenerate (graded) Lusztig algebras. See [L,C5]
and [O1], Definition 8.2. 
However the main applications
(say, the raising operators)
 require  affine ones.

Let us mention that 
equating renormalized  $\Phi'$ and $\Phi$ ([C5], Corollary 2.5 and [C2],
Appendix)
 one comes to
a Lusztig-type isomorphism [L] between 
proper completions of the algebras $\HH'$ and $\HH$.

{\it Difference-rational limit.} We follow the same procedure
(\ref\onelim), but now assume that 
$$
\eqalignno{
&X_b=q^{\la_b}, \ \la_{a+b}=\la_a+\la_b, \ \la_{[b,u]}=\la_b+u.
&\eqnu
\label\twolim\eqnum*
}
$$
It changes the result drastically.
Degenerating $\{\hT_j\}$ from (\ref\Demaz), we come to
the {\it rational  Demazure-Lusztig operators} from [C6]:
$$
\eqalignno{
&S_j\  = \   s_j^\la\ +\ 
{k_j \over \la_{\al_j}}(s_j^\la-1),\
\ 0\le j\le n,
&\eqnu
\label\hatS\eqnum*
}
$$
where by $\hw^\la$ we mean the action on $\{\la_b\}$:
$\hw^\la(\la_b)=\la_{\hw(b)}$. 
For instance,
$
S_0  =   s_0^\la + 
{k_0 \over 1-\la_\th}(s_0^\la-1).
$

We set 
$$
\eqalignno{
& S_{\hw}\equal\pi_r^\la S_{i_1}\ldots S_{i_l} \for
\hw=\pi_r s_{i_1}\ldots s_{i_l}, 
&\eqnu 
\label\Shwde\eqnum*
\cr
&\De_b\equal S_b \for b\in B,\ \De_i=\De_{b_i}.
&\eqnu
\label\Desb\eqnum*
}
$$
The definition does not depend on the particular
choice of the decomposition of $\hw\in W^b$, and
the map $\hw\mapsto S_{\hw}$ is a homomorphism. 
The operators  $\De_b$ are pairwise 
commutative.

The limit of  $\mu$ is
the {\it asymmetric Harish-Chandra 
function}:
$$
\eqalign{
&\si\ =\ \prod_{a \in R_+}
{\Ga(\la_\al+k_a)\Ga(-\la_\al+k_\al+1)
\over
\Ga(\la_\al)\Ga(-\la_\al+1)
},
}
\eqnu
\label\sigma\eqnum*
$$
where $k_\al=\ka_\al(\al,\al)/2$, $\Ga$ is
the classical $\Ga$-function.

\proclaim{ Proposition}
a) The operators $S_{\hw}$ are well-defined and  preserve
the space of polynomials $\C[\la]$ in terms of $\la_b$.
The map $\hw\mapsto S_{\hw},\ 
y_b\mapsto \la_b$ is a representation of $\HH' $.

b) The operators $S_{\hw}$ are formally unitary with respect to the inner
product 
$$
\eqalignno{
&\{f,g\}_\si\equal \int f(\la)g(\la)\si d\la,
&\eqnu
\label\siprod\eqnum*
}
$$
i.e. $\si^{-1}S_{\hw}^+\si=S_{\hw}^{-1}$ for the anti-involution $^+$
sending $\la_b\mapsto \la_b,\ \hw\mapsto \hw^{-1},$
where $\hw\in W^b$.
\label\DESI\theoremnum*
\endproclaim
{\it Proof.} Statement a) is a straightforward degeneration
of the polynomial representation of $\HH\ $.
Since $(X_b)^*=(q^{\la_b})^*=(X_b^{-1})^*$
and $q^*=q^{-1}$, we get the $*$-invariance of $\la$
and come to  (\ref\siprod).
\proofbox

Following  
Theorem \ref\faith (and [C6]),
the operators 
$$
\eqalign{
&\La_p=p(\De_{1},\ldots,\De_{n}) \for p\in 
\C[X_1^{\pm 1},\ldots, X_n^{\pm 1}]^W
}
\eqnu
\label\Laop\eqnum*
$$
are $W$-invariant 
and preserve $\C[\la]^W$ (usual symmetric polynomials in $\la$).

Let us degenerate  Proposition \ref\MACOPER (the notations
$r, m_r, \ell_r$ remain the same). The restriction of $\La_{m_r}(\De_b)$
onto $\C[\la]^W$ is as follows:
$$
\eqalign{
& \La_r= \ \sum_{w\in W/W_r}
 \prod _{\al \in \ell_r}
{\la_{w(\al)}+k_\al \over
     \la_{w(\al)}   } \, w(-b).
}
\eqnu
\label\macdeop\eqnum*
$$
In the differential case the formulas (and the corresponding
limiting procedure) are more complicated.

%
%
%
%
\section { Opdam transform } 
From now on we assume that $\R\ni k_\al>0$ for all $\al$
and keep the notation from the previous section; $i$ is the imaginary
unit, $\Re,\Im$ the real and imaginary parts.
\proclaim {Theorem}
 There exists a solution $G(x,\la)$ of the eigenvalue
problem
$$
\eqalign{
& D_b (G(x,\la))=\la_b G(x,\la),\  b\in B,\  G(0,\la)=1
}
\eqnu
\label\eigenop\eqnum*
$$
 holomorphic for all $\la$ 
and for $x$ in $\R^n+iU$ for a neighbourhood $U\subset \R^n$ of
zero. If $x\in \R^n$ then
$|G(x,\la)|\le |W|^{1/2}\exp(\max_w(w(x),\Re\la))$, so 
$G$ is bounded for $x\in \R^n$ when $\la\in i\R^n$.
The solution of (\ref\eigenop) is unique in the class of
continuously differentiable functions on $\R^n$ (for a given
$\la$).
\label\OPDAM\theoremnum*
\endproclaim

This theorem is from  [O1] (Theorem 3.15 and Proposition 6.1). 
Opdam uses that
$$
\eqalign{
& F\equal |W|^{-1}\sum_{w\in W}G(w(x),\la)
}
\eqnu
\label\fhyper\eqnum*
$$
 is a 
{\it generalized hypergeometric function}, i.e a
$W$-symmetric eigenfunction of the 
restrictions $L'_p$ of the operators $p(D_{b_1},\ldots,D_{b_n})$
to symmetric functions:
$$
L_p'F(x,\la)\ =\ p(\la_1,\ldots, \la_n)F(x,\la),
$$
where $p$ is any
$W$-invariant polynomial of $\la_i=\la_{b_i}$.
The operators $L'_p$  generalize the radial parts
of Laplace operators on the corresponding symmetric space
(see Introduction).
The normalization is the same: $F(0,\la)=1$. 
It fixes $F$ uniquely. So it is $W$-invariant with respect to $\la$
as well.

A systematic algebraic
and analytic theory of $F$-functions is due to Heck\-man and Op\-dam 
(see [HO1,H1,O2,HS]). There is 
a formula for $G$ in terms of $F$ (at least for generic $\la$)
via the operators $D_b$ from (\ref\dunk).
The positivity of $k$ implies that  
it holds for all $\la\in \C^n$.
See [O1] for a nice and simple argument (Lemma 3.14). 
This formula and  the relation to the affine
Knizhnik-Zamolodchikov equation [C5,Ma,C8,O1] are applied to 
establish the growth
estimates for $G$. Actually it  gives more than was 
formulated in the theorem
(see Corollary 6.5, [O1]).

We introduce the {\it Opdam transform}
(the first component of what he 
called ``Cherednik's transform'') as follows:
$$
\eqalignno{
&\f(f)(\la)\equal \int_{\R^n}f(x)G(-x,\la)\tau dx,
&\eqnu
\label\Optran\eqnum*
}
$$ 
for the standard measure $dx$ on $\R^n$.

\proclaim{ Proposition}
a) Let us assume that $f(x)$ are taken from the space $\C^\infty_c(\R^n)$
of $\C$-valued  compactly supported $\infty$-differentiable functions
on $\R^n$. The inner product 
$$
\eqalignno{
&\{f,f'\}_\tau\equal \int_{\R^n}f(x)f'(-x)\tau dx
&\eqnu
\label\ffprime\eqnum*
}
$$
satisfies the conditions of part b),
 Proposition \ref\DUNDEG. Namely, $\{D_b\}$ preserve $\C^\infty_c(\R^n)$
and are self-adjoint with respect to the pairing (\ref\ffprime).

b) The Opdam transforms  of such functions are analytic in $\la$
on the whole $\C^n$ and satisfy the Paley-Wiener condition. A function
$g(\la)$ is of PW-type ($g\in PW(\C^n)$) if
there exists a constant $A=A(g)>0$  such that for any $N>0$
$$
\eqalignno{
& g(\la)\le C(1+|\la|)^{-N}\exp (A|\Re\la|)
&\eqnu
\label\PWcond\eqnum*
}
$$ 
for a proper constant $C=C(N;g)$.
\label\OPPROP\theoremnum*
\endproclaim
{\it Proof.} The first claim is obvious. The Paley-Wiener condition
follows from Theorem 8.6 [O1]. The transform under consideration is actually
the first component of Opdam's transform from Definition 7.9 (ibid.).
We also omit the complex conjugation in his definition and change the
sign of $x$ (instead of $\la$). The estimates remain the same. 
\proofbox

%
%
%
%
\section { Inverse transform } 
We are coming to the inversion procedure (for positive $k$).
 The  {\it inverse transform} is well-defined for Paley-Wiener
functions $g(\la)$ on $\C^n$ by the formula
$$
\eqalignno{
&\g(g)(x)\equal \int_{i\R^n}g(\la)G(x,\la)\si d\la
&\eqnu
\label\invtran\eqnum*
}
$$ 
for the standard measure $d\la$.
The transforms of such $g$ belong to $\C^\infty_c(\R^n)$. 

 The existense readily follows from
(\ref\PWcond) and known properties
of the  Ha\-rish-Chan\-dra $c$-function (see  below).  
The embedding $\g(PW(\C^n))\subset\C^\infty_c(\R^n)$  is due to 
Opdam. It is similar to the classical
one from [He1] (see also [GV]).

Let us first discuss the shift of integration contour in (\ref\invtran).
There exists an open neighborhood
$U_+^a\subset \R^n$ of the {\it closure} $\bar{C}_+^a\in \R^n$ of the
affine Weyl chamber $C^a_+$ from (\ref\afchamber) such that
$$
\eqalignno{
&\g(g)(x)\ =\  \int_{\xi+ i\R^n}g(\la)G(x,\la)\si d\la
&\eqnu
\label\contour\eqnum*
}
$$ 
for $\xi\in U_+^a$. 
Indeed,  $k_\al>0$ and $\si$ has no singularities on $U_+^a+i\R^n$.
Then we use
the classical formulas for $|\Ga(x+iy)/\Ga(x)|$
for real $x,y, x>0$. It gives that (cf. [He1,O1])
$$
\eqalignno{
&|\si(\la)| \ \le \ C(1+|\la|)^K, \where K=2\sum_{ \al>0} k_\al,\
\la\in U^a_+ +i\R^n,
&\eqnu
\label\cestim\eqnum*
}
$$
for sufficiently big $C>0$. So the products
of PW-functions by $\si$  tend to  $0$ for $|\la|\mapsto \infty$,
and we can switch to  $\xi$.
Actually we can do this  for any integrand analytic in $ U^a_+ +i\R^n$
and approaching $0$ at $\infty$.
We come to the following:

\proclaim{ Proposition}
The conditions of part b), Proposition \ref\DESI are satisfied
for 
$$
\eqalignno{
&\{g,g'\}_\si\equal\int_{i\R^n}g(\la)g'(\la)\si d\la
&\eqnu
\label\ggprime\eqnum*
}
$$ 
in the class of PW-functions, i.e. the operators $S_{\hw}$ are
well-defined on such functions and unitary.
\endproclaim
{\it Proof.} 
It is sufficient to check the unitarity  for the generators 
$S_j= S_{s_j}\ (0\le j\le n)$ and $\pi_r^\la\ (r\in O^*)$.
 For instance, let us consider $s_0$. We will
integrate over $\xi+i\R$, assuming that  
$\xi'\equal s_\th(\xi)+\th\in U_+^a$   
and avoiding the wall $(\th,\xi)=1$.

We will apply  (\ref\siprod), the formula
$$
\int_{\xi+i\R^n}s_0(g(\la)\si) d\la \ =\
\int_{\xi'+i\R^n}g(\la)\si d\la \ =\ \int_{\xi+i\R^n}g(\la)\si d\la,
$$
and a similar formula for $(1-\la_\th)^{-1}g\si$, where $g$ is of
PW-type on $ U^a_+ +i\R^n$.
Note that $(1-\la_\th)^{-1}\si$ is regular in this domain.
One has:
$$
\eqalignno{
&\int_{\xi+i\R^n}S_{0}(g(\la))\ g'(\la)\si d\la\cr 
=\ &\int_{\xi+i\R^n}\big(s_0+k_0(1-\la_\th)^{-1}(s_0-1)\big)(g(\la))\ g'(\la)
\si d\la
\cr
=\ &\int_{\xi+i\R^n}g(\la)\ \big(s_0+k_0(s_0-1)(1-\la_\th)^{-1}\big)(g'(\la)
\si) d\la
\cr
=\ &\int_{\xi+i\R^n}g(\la)\ S_{0}^{-1}(g'(\la))\si d\la. 
&\eqnu
\label\intint\eqnum*
}
$$
Since (\ref\intint) holds for one $\xi$ it is valid for all
of them in $U_+^a$
including $0$.
The consideration of the other generators is  the same.
\proofbox

\proclaim {Main Theorem}
Given $\hw\in W^b,\ b\in B,\ f(x)\in \C^\infty_c(\R^n),\ g(\la)\in PW(\C^n),$
$$
\eqalignno{
&\hw^x(G(x,\la))\ =\ S_{\hw}^{-1}G(x,\la), \hbox{\ e.g.\ }
 X_b G\ =\ \De_b^{-1}(G),
&\eqnu
\label\gidesym\eqnum*
\cr
&\f(\hw^x(f(x)))\ =\ S_{\hw}\f(f(x)),\ 
\g(S_{\hw}(g(\la)))\ =\ \hw^x(\f(f(x))), 
&\eqnu
\label\gisym\eqnum*
\cr
&\f(X_b f(x))\ =\ \De_b(\f(f(x))),\
\g(\De_b(g(\la)))\ =\ X_b\g(g(\la)),
&\eqnu
\label\fxsym\eqnum*
\cr
&\f(D_b(f(x)))\ =\ \la_b\f(f(x)),\
\g(\la_b g(\la))\ =\ D_b(\g(g(\la))).
&\eqnu
\label\fdsym\eqnum*
}
$$
\label\GISYM\theoremnum*
\endproclaim
{\it Proof.}
Applying the intertwiners from 
(\ref\Phiprime) to $G(x,\la)$, we see that
$$
\eqalignno{
&(1 + {k_i \over \la_{\al_i} })^{-1}(s_i + {k_i \over \la_{\al_i} })(G)\ =
s_i^\la(G), \for 1\le i\le n,
&\eqnu \cr
\label\Phig\eqnum*
&(1 + {k_0 \over 1-\la_{\th} })^{-1}
(X_\th s_\th + {k_0 \over 1-\la_{\th} })(G)\ =
s_0^\la(G),\cr
&P_r'(G)\  =\  X_r\om_r^{-1}(G)=(\pi_r^\la)^{-1}(G)
 \for \ r\in O^*.
}
$$
The scalar factors on the left are necessary to preserve the
normalization $G(0,\la)=1$, so we can use the main property
of the intertwiners (\ref\Phix) and the uniqueness of $G(x,\la)$
(Theorem \ref\OPDAM). Expressing $s_j^x$
in terms of $s_j^\la$ (when applied to $G$!), we get
(\ref\gidesym) for $\hw=s_j$.   
It is obvious for $\hw=\pi_r$. Using the commutativity of $\hw^x$
and $S_{\hu}$ for $\hw,\hu\in W^b$, we establish this relation
in the general case. For $w\in W$ it
is due to Opdam.

Formula (\ref\gisym) results directly from
(\ref\gidesym) because we already  know that $\hw^x$
are unitary for $\{f,f'\}_\tau$ and $S_{\hw}$ are unitary with respect
to $\{g,g'\}_\si$ for the considered classes of functions.
See Theorems \ref\DUNDEG, \ref\DESI and (\ref\ffprime),(\ref\ggprime). 

For instance, let us check (\ref\fxsym) for $\f$, which is a particular case
of (\ref\gisym):
$$
\eqalignno{
&\f(X_bf(x))\ =\ \int_{\R^n} f(-x) X_b^{-1}G(x,\la)\tau dx \ =\cr
&\int_{\R^n} f(-x) \De_b(G(x,\la))\tau dx \ =\ 
\De_b(\f(f(x))).
&\eqnu
\label\trxb\eqnum*
}
$$
Thus the 
 Opdam transforms of the ``multiplications by the coordinates'' 
$X_b$ are the
operators $\De_b$. 

Since $D_b$ and $\la_b$ are self-adjoint for the
corresponding inner products (see Theorem \ref\DUNDEG), we
get (\ref\fdsym), which is in fact the defining property of $\f$ and $\g$. 
\proofbox

\proclaim {Corollary}
The compositions $\g\f: C_c^\infty(\R^n)\to  C_c^\infty(\R^n)$,
and $\f\g: PW(\C^n)\to  PW(\C^n)$ are multiplications by nonzero
constants. The transforms  $\f,\g$ establish isomorphisms 
between the corresponding space identifying 
$\{\ ,\ \}_\tau \and \{\ ,\ \}_\si$ up to proportionality.
\label\INVERSE\theoremnum*
\endproclaim
{\it Proof.}
The first statement readily follows from the Main Theorem.
The composition $\g\f$ sends $\C^\infty_c(\R^n)$ into itself,
is continuous (due to Opdam),
and commutes with the operators $X_b$ (multiplications by $X_b$).
So it is the multiplication by a function $u(x)$ of $\C^\infty$-type. 
It must also commute
with $D_b$. Hence $G(x,\la)u(x)$ is another solution
of the eigenvalue problem (\ref\eigenop) and $u(x)$ has to be
a constant. 

Let us  check that $\f\g$, which is
a continuous operator on $PW(\C^n)$ (for any fixed $A$)
commuting with multiplications by any $\la_b$, has to be
a multiplication by an analytic function $v(\la)$. 
Indeed, the image of $\f\g$ with respect to the
standard  Fourier transform ($k=0$)
is a continuous operator on $C_c^\infty(\R^n)$ commuting with
the derivatives $\partial/\partial x_i$. So it is the convolution
with some function and its inverse Fourier transform is the
multiplication by a certain $v(\la)$.
Since $\f\g(g)=\hbox {Const\ } g$
for any $g(\la)$ 
from $\f(C_c^\infty(\R^n))$,  $v$ is  constant.

The claim about the inner products is obvious because
 $\{\f (f) , g \}_\si$ =  $\{f ,\g(g) \}_\tau$
for $f\in C_c^\infty(\R^n)$, $g\in PW(\C^n)$. 
\proofbox

The corrolary is due to Opdam ([O1],Theorem 9.13 (1)).
He uses  Peetre's
characterization of differential operators (similar to what 
Van Den Ban and Schlichtkrull did in [BS]).
Anyway a certain nontrivial analytic argument 
is involved to check that
$\g\f$ is multiplication by a constant.

{\it The symmetric case.} The transforms can be readily
reduced to the symmetric level. If $f$ is 
$W$-invariant then 
$$
\eqalign{
&\f(f(x))\ =\  \int_{\R^n} f(x) F(-x,\la)\tau dx
}
\eqnu
\label\hartran\eqnum*
$$
for  $F$ from (\ref\fhyper).
Here we applied the $W$-sym\-met\-rization to the integrand of the
(\ref\Optran) and used that $\tau$ is $W$-invariant.
So $\f$ {\it coincides} with  the $k$-deformation of the
Harish-Chandra transform
on $W$-invariant functions 
up to a minor technical detail. The $W$-invariance of  $F$ in $\la$ 
results in the $W$-invariance of $\f(f)$.

As to $\g$, we $W$-symmetrize the integrand in the definition
with respect to $x$ and $\la$,
replacing  $G$ by  $F$ and $\si$ by
by its $W$-symmetrization. The latter
is the genuine
Harish-Chandra ``measure''
$$
\eqalign{
&\si'\ =\ \prod_{a \in R_+}
{\Ga(\la_\al+k_a)\Ga(-\la_\al+k_\al)
\over
\Ga(\la_\al)\Ga(-\la_\al)
}
}
\eqnu
\label\sigmap\eqnum*
$$
up to a coefficient of proportionality.

Finally, given a $W$-invariant function $f\in\C^\infty_c(\R^n),$
$$
\eqalignno{
& f(x)\ =\ \hbox{Const} \int_{i\R^n}g(\la)F(x,\la)\si' d\la \for g=\f(f).
&\eqnu
\label\harish\eqnum*
}
$$
A similar formula holds for $\g$.
See [HC,He2] and [GV], Ch.6 for the classical theory.

As an application, we are able to
 calculate the Fourier transforms of the operators
$p(X)\in \C[X_b]^W$ (symmetric Laurent polynomials acting by multiplication) 
in the  Harish-Chandra theory and its $k$-deformation. 
They are exactly the 
operators $\La_p$ from Section 4.
In the minuscule case, we get  formulas (\ref\macdeop).
We mention that in  [O1] and other papers 
the pairings serving the Fourier transforms are hermitian.
Complex conjugations can be added to ours. 

We hope that the method used in this paper can be generaized to
negative $k$ and to other classes of functions (cf. [BS,HO2]). 
Relations (\ref\gidesym) considered as difference equations for $G(x,\la)$
with respect to $\la$ may help with the growth estimates via 
the theory of difference equations and the equivalence with
difference Knizhnik-Zamolodchikov equations.
However  I believe
that the main progress here will be connected with 
the difference Fourier transform.

%
%
%
%
%
\AuthorRefNames [BGG]
\references
\vfil

[BS]
\name{E.P. van den Ban}, and \name{H. Schlichtkrull},
{The most continuous part of the Plancherel decomposition for
a reductive symmetric space}, Preprint (1993).

[B]
\name{N. Bourbaki},
{ Groupes et alg\`ebres de Lie}, Ch. {\bf 4--6},
Hermann, Paris (1969).

[C1]
\name{I. Cherednik},
{Difference Macdonald-Mehta conjectures},
IMRN (1997).

[C2]
\bibline, 
{ Intertwining operators of double affine Hecke algebras},
Selecta Math. (1997).

[C3]
\bibline, 
{ Macdonald's evaluation conjectures and
difference Fourier transform},
Inventiones Math. {122} (1995),119--145.

[C4]
\bibline, 
{ Nonsymmetric Macdonald polynomials },
IMRN {10} (1995), 483--515.

[C5]
\bibline, 
{A unification of Knizhnik--Zamolodchikov
and Dunkl operators via affine Hecke algebras},
Inventiones Math. {  106}:2  (1991), 411--432.

[C6]
\bibline,
{ Quantum Knizhnik-Zamolodchikov equations and affine
root systems},
Commun. Math. Phys. { 150} (1992), 109--136.

[C7]
\bibline,
{Elliptic quantum many-body problem and double affine
 Knizhnik - Zamolodchikov equation},
Commun. Math. Phys. {169}:2  (1995), 441--461.

[C8]
\bibline,
{ Integration of quantum many- body problems by affine
Knizhnik--Za\-mo\-lod\-chi\-kov equations}, 
Advances in Math. {106}:1 (1994), 65--95
(Pre\-print RIMS--{  776} (1991)).

[C9]
\bibline,
{ Double affine Hecke algebras, 
Knizhnik- Za\-mo\-lod\-chi\-kov equa\-tions, and Mac\-do\-nald's 
ope\-ra\-tors},
IMRN (Duke M.J.) {  9} (1992), 171--180.

[C10]
\bibline,
{Lectures on affine Knizhnik-Zamolodchikov
equations, quantum many-body problems, Hecke
algebras, and Macdonald theory}, Preprint RIMS-1144,\  in
MSJ Memoirs
(1997).

[C11]
\bibline,{ Double affine Hecke algebras and  Macdonald's
conjectures},
Annals of Mathematics {141} (1995), 191--216.

[D]
\name{C.F. Dunkl},
{ Hankel transforms associated to finite reflection groups},
Contemp. Math. {138} (1992), 123--138.

[GV]
\name{R. Gandolli}, and \name {V.S. Varadarajan},
{Harmonic analysis of spherical functions on
real reductive groups},
Springer-Verlag (1988).

[HC]
\name {Harish-Chandra}
{Discrete series for semisimple Lie groups},
Acta Math. {116} (1986), 1-111.

[H1]
\name{G.J. Heckman},
{ Root systems and hypergeometric functions II},
Comp. Math. {64} (1987), 353-374.

[HO1]
\name{G.J. Heckman}, and \name{E.M. Opdam},
{ Root systems and hypergeometric functions I},
Comp. Math. {64} (1987), 329--352.

[HO2]
\bibline,
{Harmonic analysis for affine Hecke algebras},
Preprint (1996).

[HS]
\name{G.J. Heckman}, and \name {H. Schlichtkrull},
{Harmonic analysis and special functions on symmetric spaces},
Persp. in Math., Acad. Press (1994).

[He1]
\name {S. Helgason},
{Groups and geometric analysis},
Academic Press, New York (1984).

[He2]
\name {S. Helgason},
{An analogue of the Paley-Wiener theorem for the Fourier
transform on certain symmetric spaces},
Math. Ann., {165} (1966), 297-308.

[J]
\name{M.F.E. de Jeu},
{The Dunkl transform }, Invent. Math. {113} (1993), 147--162.

[Kn]
\name {F. Knop},
{Integrality of two variable Kostka functions},
Preprint (1996). 

[KS]
\name {F. Knop}, and \name{S. Sahi},
{A recursion  and a combinatorial formula for Jack
polynomials},
Preprint (1996).

[L]
\name {G. Lusztig},{ Affine Hecke algebras and their graded version},
J. of the AMS { 2}:3 (1989), 599--685.

[M1]
\name{I.G. Macdonald}, {  A new class of symmetric functions },
Publ.I.R.M.A., Strassbourg, Actes 20-e Seminaire Lotharingen,
(1988), 131--171 .

[M2]
\bibline, {  Orthogonal polynomials associated with root 
systems},Preprint(1988).

[Ma]
\name {A. Matsuo}, {Knizhnik-Zamolodchikov type equations and zonal
spherical functions},
Inv.Math. {110}, (1992), 95--121.

[O1]
\name{E.M. Opdam}, 
{Harmonic analysis for certain representations of
graded Hecke algebras}, 
Acta Math. {175}, (1995), 75--121.

[O2]
\bibline, 
{  Some applications of hypergeometric shift
operators}, Inventiones Math. {  98} (1989), 1--18.

[R]
\name{S.N.M. Ruijsenaars},
{Complete integrability of relativistic Calogero-Moser
systems and elliptic function identities}, Commun. Math. Phys.
{ 110} (1987), 191--213.

\endreferences
\bye

%